\documentclass [12pt] {article}

\usepackage[cp850]{inputenc}
\usepackage{amsmath,amsfonts,colortbl}

\usepackage{xspace,ifthen,a4,a4wide}
\usepackage{Umlaute,Farbe,German,Schrift}

\usepackage[engl,aw]{F}

\def \zitC {Carter [1972]\xspace}

\def \zitA {Conway \& al. [1984]\xspace}

\def \zitG {Gr\"uninger [2007]\xspace}

\def \zitHS {Havas \& Sims [1999]\xspace}

\def \zitH {Humphreys [1975]\xspace}

\def \zitJW {Jansen \& Wilson [1996]\xspace}

\def \zitK {Kantor [1981]\xspace}

\def \zitL {Lyons [1972]\xspace}

\def \zitMN {Meyer \& Neutsch [1984]\xspace}

\def \zitMNP {Meyer \& al. [1985]\xspace}

\def \zitMo {Moore [1897]\xspace}

\def \zitNM {Neutsch \& Meyer [1989]\xspace}

\def \zitS {Sims [1973]\xspace}

\def \zitSt {Steinberg [1962]\xspace}

\def \zitWoth {Woldar [1984]\xspace}

\def \zitWop {Woldar [1987]\xspace}

\usepackage {Umlaute}
\usepackage {graphicx}
\usepackage {epsfig}
\usepackage {amssymb}

\begin {document}

\begin {center}
{\bf \LARGE \color {blue} Lyons Taming} \\[5 mm]
{\large \it by} \\[5 mm]
{\bf \LARGE \color {blue} Wolfram Neutsch} \\[5 mm]
{\small \bf \color {red} University of Bonn [retired]} \\[5 mm]
{\color {black} \large \it Vorster Str. 29 \\
\it D--41169 M\"onchengladbach \\
\it Deutschland/Germany} \\[5 mm]
{\sf \small \bf \color {red} Wolfram\_Neutsch@gmx.de} \\[4 mm]
{\color {black} \small \bf 2010 Mathematics Subject Classification: 20D08, 20C34, 20C20, 20F05} \\
{\color {black} \small \bf Key words: Sporadic Lyons group, 5-representation, Kantor geometry} \\[5 mm]
{\sf \small \bf \color {red} 16. August 2023} \\[1 cm]
\end {center}

Based on Kantor's geometry, we give a new highly symmetric construction of Lyons' sporadic group $Ly$ via its minimal representation over $\F_5^{111}$, thus obtaining elementary existence proofs
for both the group and the representation at one stroke.

\thispagestyle {empty}

\section {Introduction: History of the Lyons group}

\label {a.hist}

As part of a more general classification problem, \zitL studied the question whether there exists a simple group $\Gamma$ with an involution $z \in \Gamma$ such that the centralizer $C_\Gamma(z)$ is isomorphic to the Schur double cover $2 \Schur A_{11}$ of the alternating group $A_{11}$ on 11 letters (for the sake of clarity, we here and elsewhere modify the notation in order to achieve a more homogeneous nomenclature).

Lyons showed that such a putative $\Gamma$ must have several very precisely describable properties. For example, the local subgroup structure and the order of $\Gamma$ are uniquely determined, and Lyons even succeeded in constructing the complete character table of $\Gamma$. Moreover, he was able to demonstrate that the largest proper subgroups of $\Gamma$ form a single conjugacy class and are isomorphic to the simple Chevalley group $G_2(5)$.

This extremely important assertion was taken up by Sims who chose a set of four elements in $G_2(5)$ and deduced a number of relations in them which generate $G_2(5)$ as an abstract group. Then he showed that every group like $\Gamma$ must contain a fifth element fulfilling (together with the four generators already given) 14 additional conditions. Furthermore, all relations combined form a presentation for $\Gamma$.

These results simultaneously establish existence and uniqueness (up to isomorphy) of $\Gamma$, which we shall call from now on in honor of its discoverer the {\bf Lyons group}, conventionally described by the symbol $Ly$.

Unfortunately, \zitS decided to write down only the 14 additional relations explicitly, while he defined the $G_2(5)$-generators as matrices in a faithful 7-dimensional representation of $G_2(5)$ over $\F_5$. The complete set of relations required to characterize $Ly$ was published only much later in \zitHS. This undue delay caused many inconveniences in the meantime, as we shall see shortly.

A major breakthrough in our understanding of the fundamental properties of $Ly$ was achieved in the seminal paper of \zitK, in which the geometric properties of the Lyons group in characteristic 5 are elucidated. It turned out that many features known from Chevalley theory have direct parallels in $Ly$.

In short terms, Kantor introduced the ground field $\F_5$ and consequently defined as {\bf Borel groups} the Sylow-5-normalizers in $Ly$. A {\bf parabolic group} then is a supergroup of some Borel group with the single exception of $Ly$ itself.

The {\bf objects} of {\bf Kantor's geometry} are the maximal parabolics which fall into three conjugacy classes under the action of $Ly$. Kantor divided his objects accordingly into three families which may be distinguished by their isomorphism types: "{\bf points}" are isomorphic to $G_2(5)$, "{\bf lines}" to $5^{1+4}:4S_6$ and "{\bf planes}" to $5^3 \cdot SL_3(5)$ (in {\sc Atlas} notation, cf. \zitA).

The geometry itself consists of these objects and a relation called {\bf incidence}. A set $M$ of objects is incident (or a {\bf flag}) if their intersection contains a Borel subgroup. Note that the empty set conventionally has the universal group $Ly$ as intersection; hence $\leer$ by definition is a flag. All nonempty inciding sets have parabolic groups as intersections.

Every maximal flag contains precisely one object of each type, and all maximal flags are conjugate in the Lyons group. The same holds, by the way, in the above-mentioned Chevalley groups. The stabilizer of a maximal flag is equal to its intersection, namely the unique Borel group it contains.

The Borel groups are -- just like those of $A_2(5)$ or $G_2(5)$ -- split extensions of the underlying Sylow-5-group $S$, say, with 2-groups isomorphic to the direct product of two copies of the multiplicative group $\F_5^\times \iso 4$.

Any complement of $S$ in its normalizer is called a {\bf torus} of $Ly$. This is also in perfect analogy with the Chevalley geometries. As is easily seen, all tori are conjugate.

Of utmost importance for Kantor's theory, however, is another concept borrowed from Chevalley geometry. The {\bf apartment} $\frak A = \frak A(T)$ associated with some particular torus $T$ consists of the exactly 12 points, 36 lines, 24 planes and 144 flags stabilized by $T$.

The $\frak A$-objects are consequently permuted by the {\bf torus normalizer} $N = \n_{Ly}(T)$. This leads to a natural permutation action of the {\bf Weyl group} $W = N/T$ (another term familiar from Chevalley theory) on the apartment which is transitive on each type of objects and even sharply transitive on the maximal flags.

Kantor gave each of the 12 points in $\frak A$ a {\bf name} composed of two {\bf coordinates}, a number taken from the set $\{1,2,3,4\}$ and a letter from $\{a,b,c\}$. Each line connects two points which differ in both coordinates; each plane incides with exactly three such points and the three lines formed by them. The Weyl group is $S_{\{1,2,3,4\}} \times S_{\{a,b,c\}}$.

\begin {center}

\abb {16cm} {16cm} {apart} {\bf Apartment (\zitK)} \label {abb.k}

\end {center}

The incidence relation endows $\frak A(T)$ with the structure of a 2-dimensional simplicial complex whose 0-, 1- and 2-simplices are the points, lines and planes, respectively. A picture of $\frak A(T)$ is displayed in Fig. \ref {abb.k}. Of course, points carrying the same names have to be identified.

In general, the description of a finite group by generators and relations is not very useful for calculating inside the group. This is particularly true for the Sims presentation of $Ly$. It is therefore of great interest to find a faithful representation of the Lyons group. Alas, it is not easy to construct a suitable one.

The simplest nontrivial permutation representations of $Ly$ are those by multiplication on the cosets (or equivalently conjugation on the conjugates) of the two largest proper subgroups, namely $G_2(5)$ and the normalizer of some $3a$-element which is of isomorphism type $(3 \Schur Mc):2$, where $Mc$ denotes the sporadic group discovered by McLaughlin.

Both, however, are rather complicated, not so much due to their huge degrees (8835156 and 9606125, respectively), but because of their large rank (5 in both cases). All other permutation representations are even more hopeless.

Ordinary matrix representations are also quite formidable. As the character table reveals, the minimal degree is 2480, attained for instance by representations over $\C$ or $\Q(\sqrt {-11}$) and twice as large for matrices over the real or rational numbers. Therefore the attention soon shifted to modular representations in some prime characteristic $p$ dividing the group order. namely $p \in \{ 2,3,5,7,11,31,37,67 \}$.

\zitMN launched their investigation into the relevant problem by attempting to establish lower bounds for the degree $n$ of some nontrivial $p$-representation of $Ly$ (which due to the simplicity of $Ly$ will automatically be faithful).

So let $D$ be an injective homomorphism from $Ly$ into the general linear group $GL_n(\K)$ over some field $\K$ of characteristic $p$ with $n$ as small as possible and $\phi$ the corresponding Brauer character. Since $n$ is minimal, $D$ (and $\phi$) are irreducible, and by a standard result of Brauer theory $D$ is (up to similarity) uniquely determined by $\phi$.

According to its definition, $\phi$ is a complex class function on the set of all $p$-regular elements of $Ly$, and the restriction to any (proper) subgroup $U < Ly$ must yield a valid Brauer $p$-character. Executing this idea for $U$ isomorphic to the (maximal) subgroups $G_2(5)$, $(3 \Schur Mc):2$, $2 \Schur A_{11}$ and $3^5:(2 \times M_{11})$, where $M_{11}$ is the smallest sporadic Mathieu group, led to severe restrictions imposed on $\phi$.

Meyer and Neutsch deduced from their calculations that $n$ must be at least equal to 111 if $p = 5$, while in any other characteristic the lower limit is 124. As we now know, the latter value is a vast underestimate; the smallest degree faithful matrix representation of the Lyons group in any characteristic other than 5 is 651-dimensional over $\F_3$ and was first constructed by \zitJW.

On the other hand, it was impossible to rule out an injective degree-111-representation for $p = 5$. In contrast, the existence of such a $D$ and its Brauer character was highly probable since the calculations led to exactly three possible solutions, differing only on the conjugacy classes $67a, 67b, 67c$ of order 67 elements in $Ly$.

Two of the putative characters seemed implausible as they would force the fixed spaces of the elements of order 67 to be very large. Therefore the authors conjectured that the third candidate should be the unique irreducible Brauer-character of the Lyons group with degree 111. This amounts to stating that the absolute minimal representation of $Ly$ is unique and given by 111-dimensional matrices over some field of characteristic 5.

At this point it is appropriate to stress that \zitWoth independently calculated most values of $\phi$ making use of a delicate local analysis. In this way he came to the same conclusion concerning the absolute minimal representation of $Ly$.

The truth of the above-mentioned conjecture was established one year later by \zitMNP who constructed the proposed representation $D$ explicitly and in that manner demonstrated existence and uniqueness of $D$ (hence of $\phi$ as well). Since the proof is highly technical and the details are not important for our present purpose, we shall content ourselves with briefly sketching the main ideas of the investigation.

First, Parker found two generators $a$ and $b$ in $SL_{111}(5)$ spanning a group isomorphic to $2 \Schur A_{11}$ with the correct Brauer character. Then he obtained a subgroup $M \iso M_{11}$ in $\la a,b \ra$ and determined all involutions $c \in SL_{111}(5)$ centralizing $M$ and fulfilling the additional condition $\la \z( \la a,b \ra), c \ra \iso S_3$.

There are only a palmful solutions and it easy to eliminate all $c$ which are inadmissible because $\la a,b,c \ra$ would contain elements with orders not occurring in the Lyons group. It turned out that all remaining $c$ lead to equivalent groups.

Hence the desired representation -- if it exists at all -- is unique up to similarity. At this point the values of $\phi$ can be completed, at the same time showing that the conjectured Brauer character is the only possible. Moreover the existence problem for $D$ reduces to the statement $\la a,b,c \ra \iso Ly$.

In order to prove this assertion, Meyer and Neutsch constructed a subgroup $G \iso G_2(5)$ in $\la a,b,c \ra$. This part of the proof was rather cumbersome because it required some guesswork and several properties of $G_2(5)$ as a Chevalley group, among them the structure of the (minimal) 7-dimensional 5-representation as full automorphism group of the Graves-Cayley-Dickson algebra over $\F_5$, an explicit system of root subgroups and the Steinberg presentation. Quite a number of base transformations were required to complete this step. Thereafter in a faithful 7-dimensional of $G$ four elements corresponding to Sims' first generators were identified. The remaining task, namely constructing the fifth Sims generator and verifying that the 14 proper Sims relations hold, is then easy.

Concerning a detailed account of the proof, in particular the numerous intermediate calculations (mainly base changes), the interested reader is referred to the original paper.

The degree-111-representation $D$ of the Lyons group over $\F_5$ has since found a great many of applications, for instance to the classification of the maximal subgroups of $Ly$. This was achieved by Wilson [1984,1985]. With purely group theoretical means, the maximal local groups had already been found before in Woldar's above-mentioned dissertation, cf. also the summary in \zitWop.

As we saw above, \zitK carried over most concepts familiar from Chevalley geometries to the sporadic Lyons group. Nevertheless, two of the most important ones withstood his efforts, namely root
(vector) systems and root groups. In order to close these gaps, \zitNM first modified Kantor's description of his objects (henceforth distinguished by attaching a subscript "K" for "Kantor") and
replaced them by smaller groups if possible.

In fact, the structure of the 5-geometry depends only on the way in which $Ly$ acts via conjugation on the objects, that means on their normalizers. Any object $O_K$ (i.\,e. a point $P_K$, a line $L_K$ or a plane $F_K$) in the sense of Kantor is maximal in the simple group $Ly$ and hence self-normalizing: $\n_{Ly}(O_K) = O_K$. We may therefore substitute $O_K$ by another group with normalizer $O_K$. Suited for that purpose are just the $O$ with $1 < O \nteq O_K$, and it is advisable to choose $O$ as small as allowed by this condition.

The authors follow this prescription which leads to considerable simplifications for planes and lines, but not for points. In the three relevant cases one finds:

Since $P_K \iso G_2(5)$ is simple, there is no alternative, and one has to be content with $P = P_K$. On the other hand, a Kantor line $L_K \iso 5^{1+4}:4S_6$ possesses many nontrivial normal subgroups, the smallest being $L = Z(O_5(L_K)) = Z([L_K,L_K]) \iso 5$. Similarly, for a plane one finds $F_K \iso 5^3 \cdot SL_3(5)$ and thus gets $F = O_5(F_K) \iso 5^3$. In the latter two cases, the new objects are $5a$-pure (and elementary-abelian).

We once and for all fix a torus $T$ and without further comment restrict from now on our analysis to the objects (according to the new definitions), which belong to the apartment $\frak A = \frak A(T)$. The sets of points, lines, and planes will be denoted by the symbols $\frak P$, $\frak L$ and $\frak F$ (in that order).

\zitNM then find that there seems to be no meaningful interpretation of the notion "root (vector) system" in $Ly$, at least not as a finite set of vectors with the usual symmetry properties in some Euclidean space. This, however, is merely a minor difficulty, because root systems in Chevalley groups are mainly (if not exclusively) used for bookkeeping, i.e. as a means to classify or enumerate the root subgroups.

This lack is more than outweighed by the construction of {\bf root groups} in $Ly$. The main idea behind their construction stems from the observation that tori in $A_2(5)$, $G_2(5)$ and the Lyons group itself may be viewed as those subgroups which are writable as direct products of $r$ copies of the multiplicative group of the underlying field, where the {\bf rank} $r$ is maximum possible. In the three cases at hand, one gets $r = 2$, and, consequently, all tori are of the form $(F_5^\times)^r \iso 4^2$. Vice versa, since all $4^2$-subgroups are conjugate, each of them may be considered as a torus.

The same procedure can be applied to the root group. Chevalley theory yields that they are isomorphic to the additive group of the ground field, here $\F_5^+ \iso 5$, and normalized by $T$. Carrying this over to $Ly$, the root groups turn out as identical with the 36 lines in the apartment $\frak A$.

This highly satisfying result is amended by the introduction of a system of {\bf roots} for the Lyons group, namely a certain set of elements, one generator for each root group. The authors define them in the most symmetric way, to wit, as an orbit under conjugation by a certain group acting simply transitively on the 36 lines.

The 12 root groups corresponding to the lines in the hexagonal vicinity of some point $P \in \frak P$ generate $P \iso G_2(5)$ as a group. This necessitates to introduce a further new concept. The set $\frak L$ and {\it a fortiori} the apartment $\frak A$ as a simplicial complex are endowed with an {\bf orientation} as follows:

Each line $L \in \frak L$ connects two points $P,Q \in \frak P$ whose number and letter coordinates are different. In particular, exactly two of the three letters $a,b,c$ occur in the name of $L$. Depending on them, we endow $L$ with a {\bf direction} from $a$ to $b$, from $b$ to $c$ or from $c$ to $a$, depicted in Fig. \ref {abb.nm} by the arrows. In this manner we may denote the lines uniquely as ordered pairs like $L = (P,Q)$, where the orientation is from $P$ to $Q$.

In $P$, a certain set of relations holds, forming a Steinberg presentation of $G_2(5)$. Moreover, the line set $\frak L$ naturally splits into three {\bf parallel systems} of 12 members each, as indicated by the colouring in Fig. \ref {abb.nm}. The group spanned by any parallel system is the full centralizer of one of the three torus involutions and therefore isomorphic to $2 \Schur A_{11}$.

\begin {center}

\abb {16cm} {16cm} {pfeil} {\bf Apartment with orientation (\zitNM)} \label {abb.nm}

\end {center}

Combining the 12 Steinberg presentations for the points with 3 suitable sets of relations defining the centralizers as groups generated by the parallel systems yields a presentation for some covering group of $Ly$, since by construction all relations evidently hold in $Ly$.

\zitNM conjectured that this group in fact be the Lyons group itself, but they were unable to prove this.

In spite of several attempts by a number of scholars (and methods), the problem remained open for as long as 18 years.

It was finally settled by \zitG who demonstrated that the Lyons group possesses even a slightly stronger presentation (with the same generators but fewer relations). His successful approach consists of constructing $Ly$ as an amalgam of the subgroups visible in the given data.

The new presentation has definite advantages as compared to Sims' classical result: It is geometric in spirit while the Sims generators and relations are constructed {\it ad hoc} in order to get a manageable set of equations. Much more important is, however, that the amalgam method can be carried through entirely by hand; the original proof by \zitHS required heavy computer calculations. The larger number of generators and relations of the new Ansatz is not a major drawback since it is more than outweighed by the high symmetry.

\section {Basic definitions}

\label {a.basdef}

As already stated, we throughout apply (a slightly modified form of) the almost universally accepted {\sc Atlas} nomenclature system introduced in \zitA, especially for groups and conjugacy classes. In particular, a Chevalley group is named by its underlying Lie algebra and the order of the ground field like $A_2(5)$ or $G_2(5)$.

Some types of group extensions are: $A.B$ is a general extension of $A$ by $B$, i.\,e. a group containing a normal subgroup equal (or isomorphic) to $A$ with quotient group isomorphic to $B$. More precisely, a split resp. nonsplit extension will be denoted as $A:B$ resp. $A \cdot B$, and the direct product is $A \times B$.

Finally, the symbol $A \Schur B$ always means a Schur extension of $A$ by $B$ (or a Schur $A$-cover of $B$), namely a group extension in which the normal subgroup $A$ is contained in the center $Z(A \Schur B)$ and in the commutator group $[A \Schur B, A \Schur B]$ of the whole group.

The symmetric and alternating groups on some finite set $M$ are $S_M$ and $A_M$; in case $M = \{ 1, \ldots, n \}$ with $n \in \N$ we simply write $S_n$ or $A_n$.

Recall furthermore that Kantor objects (points, lines, planes) are always restricted to the apartment $\frak A = \frak A(T)$ of some torus $T$ chosen in advance. The same holds for the incidence relation and for flags.

The sets of all 12 points, 36 lines and 24 planes belonging to $\frak A(T)$ are called $\frak P$, $\frak L$ and $\frak F$ in that order. All lines (and planes) are considered to be oriented in the above-mentioned way; cf. Fig. \ref  {abb.nm}.

We often write a line $L \in \frak L$ in the form $L = (P,Q)$. Then $P,Q \in \frak P$ are the two points incident with $L$, and the line is directed from $P$ to $Q$.

All vectors and matrices are constructed over the ground field $\F_5 = \mng {0,1,2,3,4}$. The {\bf representation space} of the group $X$ to be considered later-on is $V = \F_5^{111}$; hence the matrices operating on $V$ are elements of $\F_5^{(111,111)}$. We endow $V$ with the {\bf canonical basis} $\menge {e_i} {1 \le i \le 111}$, where $e_i$ is the $i^{th}$ row of the 111-dimensional unit matrix over $\F_5$.

Furthermore, we partition the {\bf index set} $\mng {1, \ldots, 111}$ into 16 {\bf sections} $\sigma_1, \ldots \sigma_{16}$ of consecutive numbers whose parameters are given in Table \ref {t.dim}. For each $I \in \mng {1, \ldots, 16}$, the $I^{th}$ section is $\sigma_I = (\alpha_I, \alpha_I+1, \ldots, \omega_I-1, \omega_I)$ and in particular has {\bf length} $d_I = \omega_I - \alpha_I + 1$.

\begin {center}

\begin {tabelle}

\label {t.dim}

\begin {tabular} {||c||c|c|c|c|c|c|c|c|c|c|c|c|c|c|c|c||} \hline
$I$ \rule {0mm} {4.4mm} & 1 & 2 & 3 & 4 & 5 & 6 & 7 & 8 & 9 & 10 & 11 & 12 & 13 & 14 & 15 & 16 \\ \hline
$\alpha_I$ \rule {0mm} {4.4mm} & 1 & 10 & 16 & 22 & 28 & 35 & 42 & 49 & 56 & 63 & 70 & 77 & 84 & 91 & 98 & 105 \\ \hline
$\omega_I$ \rule {0mm} {4.4mm} & 9 & 15 & 21 & 27 & 34 & 41 & 48 & 55 & 62 & 69 & 76 & 83 & 90 & 97 & 104 & 111 \\ \hline
$d_I$ \rule {0mm} {4.4mm} & 9 & 6 & 6 & 6 & 7 & 7 & 7 & 7 & 7 & 7 & 7 & 7 & 7 & 7 & 7 & 7 \\ \hline
\end {tabular}

\end {tabelle}

\end {center}

For $1 \le I \le 16$, we define $E_I$ as the subspace of $V$ spanned by $\menge {e_i} {i \in \sigma_I}$. Note that $V$ is the direct vector space sum

\begin {equation}
V = E_1 \oplus E_2 \oplus \ldots \oplus E_{15} \oplus E_{16}
\end {equation}

and $\dim E_I = d_I$ for all $I$.

This leads to a decomposition of all vectors $v \in V$ into 16 {\bf parts}:

\begin {equation}
v = (v_{[1]}, \ldots, v_{[16]}) = (v_1, \ldots, v_9 | \quad ... \quad | v_{105}, \ldots, v_{111})
\end {equation}

The $(i,j)$-{\bf component} of the matrix $x$ is denoted by the standard terminology as $x_{ij}$, while $x_{[IJ]}$ means the $(I,J)$-{\bf block} of $x$, i.\,e. the submatrix

\begin {equation}
x_{[IJ]} = \left|
\begin {matrix}
x_{\alpha_I, \alpha_J} & x_{\alpha_I, \alpha_J+1} & \ldots & x_{\alpha_I, \omega_J-1} & x_{\alpha_I, \omega_J} \\
x_{\alpha_I+1, \alpha_J} & x_{\alpha_I+1, \alpha_J+1} & \ldots & x_{\alpha_I+1, \omega_J-1} & x_{\alpha_I+1, \omega_J} \\
\ldots & \ldots & \ldots & \ldots & \ldots \\
x_{\omega_I-1, \alpha_J} & x_{\omega_I-1, \alpha_J+1} & \ldots & x_{\omega_I-1, \omega_J-1} & x_{\omega_I-1, \omega_J} \\
x_{\omega_I, \alpha_J} & x_{\omega_I, \alpha_J+1} & \ldots & x_{\omega_I, \omega_J-1} & x_{\omega_I, \omega_J}
\end {matrix}
\right|
\end {equation}

In this manner, all matrices $x \in \F_5^{(111,111)}$ are decomposed into $16 \times 16$ blocks $x_{[IJ]}$, $I,J \in \mng {1, \ldots, 16}$:

\begin {equation}
x = \left| \text{\begin {tabular} {c|c|c|c|c}
$x_{[1,1]}$ & $x_{[1,2]}$ & \ldots & $x_{[1,15]}$ & $x_{[1,16]}$ \\ \hline
$x_{[2,1]}$ & $x_{[2,2]}$ & \ldots & $x_{[2,15]}$ & $x_{[2,16]}$ \\ \hline
\ldots & \ldots & \ldots & \ldots & \ldots \\ \hline
$x_{[15,1]}$ & $x_{[15,2]}$ & \ldots & $x_{[15,15]}$ & $x_{[15,16]}$ \\ \hline
$x_{[16,1]}$ & $x_{[16,2]}$ & \ldots & $x_{[16,15]}$ & $x_{[16,16]}$
\end {tabular}} \right|
\end {equation}

Permutations of the index set $\{1, \ldots, 111\}$ are described in the usual manner as cycle products where each cycle is enclosed in round brackets.

A particular subclass are the {\bf block permutations}. These consist in {\it en bloc} permutations of the 16 parts of the above partition and are denoted by the analogously written permutation of the section numbers. To avoid misinterpretations, in this case we enclose the cycles in square brackets. Of course, this construction is restricted to those permutations (henceforth called {\bf admissible}) where each section has the same length as its image.

Hence an element of $S_{16}$ is admissible if and only if it respects the three subsets $\{1\}$, $\{2,3,4\}$ and $\{5, \ldots, 16\}$. These block permutations constitute the subgroup

\begin {equation}
S_{(1|3|12)} = S_{\{1\}} \times S_{\{2,3,4\}} \times S_{\{5,6,7,8,9,10,11,12,13,14,15,16\}} \iso S_1 \times S_3 \times S_{12}
\end {equation}

of $S_{16}$. \\

{\bf Example}

The block permutation $[1,2]$ is not defined since $d_1 = 9 \ne 6 = d_2$. On the other hand, $[2,3]$ is admissible and represents the blockwise interchange of $\sigma_2 = \{10, \ldots, 15\}$ with $\sigma_3 = \{16, \ldots, 21\}$. Written out in length, it therefore has the form

\begin {equation}
[2,3] = (10,16)(11,17)(12,18)(13,19)(14,20)(15,21)
\end {equation}

\ \\The {\bf unit matrix} in $\F_5^{(m,m)}$ will be denoted by $1_m$; if the dimension $m$ is obvious from the context we often simply write $1$. {\bf Zero vectors} and {\bf zero matrices} are always called $0$, regardless of the dimensions.

For matrices $A$ and $B$ of any format we describe by $A \oplus B$ their {\bf direct sum} and by $A \otimes B$ the {\bf Kronecker product}.

Let $I,J \in \{ 1, \ldots, 16 \}$. A matrix $B \in \F_5^{(111,111)}$ is {\bf block monomial} if it contains exactly one nonzero block in every block line and every block row. In other words, $B$ is block monomial if and only if for each $I$ there is a unique $J$ and for each $J$ a unique $I$ such that $B_{[IJ]} \ne 0$.

A block permutation matrix is block monomial with all nonzero blocks equal to the unit matrix of the appropriate dimension, and vice versa. $B$ is {\bf block diagonal} if all off-diagonal blocks vanish: $B_{[IJ]} = 0$ for $I \ne J$. Note that some or all diagonal blocks may be $0$. We occasionally use the abbreviation

\begin {equation}
\op {diag} (a_1, \ldots, a_m) = \bigoplus_{i=1}^m (a_i) = (a_1) \oplus \ldots \oplus (a_m)
\end {equation}

for the {\bf diagonal matrix} with diagonal entries $a_i$ as given and similarly

\begin {equation}
\op {Diag} (B_1, \ldots, B_{16}) = \bigoplus_{I=1}^{16} B_I = B_1 \oplus \ldots \oplus B_{16}
\end {equation}

for the {\bf block diagonal matrix} with diagonal blocks $B_I$.

Finally, $B$ is {\bf block scalar} if it is block diagonal and all diagonal blocks are scalar matrices. Thus the block scalar matrices are those of the form

\begin {equation}
\op {Diag} (c_1 1_{d(1)}, \ldots, c_{16} 1_{d(16)}) = \bigoplus_{I=1}^{16} \left( c_I \cdot 1_{d(I)} \right) = (c_1 \cdot 1_{d(1)}) \oplus \ldots \oplus (c_{16} \cdot 1_{d(16)})
\end {equation}

with coefficients $c_I \in \F_5$.

Trivially, every block monomial matrix $x$ possesses a unique decomposition $x = x_D x_P$ into a block diagonal matrix $x_D$ and a block permutation matrix $x_P$.

As is well-known, the transition from a complex Lie algebra to the associated Lie group is effected with the help of the exponential mapping which can be evaluated by the uniformly converging power series

\begin {equation}
\exp u = \Sum_{n=0}^\infty \frac {u^n} {n!} = 1 + u + \tfrac 1 2 u^2 + \tfrac 1 6 u^3 + \tfrac 1 {24} u^4 + \tfrac 1 {120} u^5 + \ldots
\end {equation}

where $u$ is some matrix over $\C$.

In quite a similar way, the canonical generators of a Chevalley group are unipotent elements of the form $\exp u$ with certain nilpotent matrices $u$ over some field $\K$ of prime characteristic $p$, say.

Under these circumstances the above development has to be modified slightly since, beginning with the term proportional to $u^p$, all coefficients in the above series have denominators divisible by $p$ and hence over $\K$ are clearly undefined. Therefore we must restrict the application of the exponential function to arguments $u$ fulfilling the additional condition $u^p = 0$. This allows to truncate the series to the finite (polynomial) expression

\begin {equation}
\exp_p u = \Sum_{n=0}^{p-1} \frac {u^n} {n!} = 1 + u + \ldots + \tfrac 1 {(p-2)!} u^{p-2} + \tfrac 1 {(p-1)!} u^{p-1} \qquad \qquad \text {(if $u^n = 0$)}
\end {equation}

Just the same phenomenon occurs for the inverse function, namely the logarithm, given by

\begin {equation}
\log (1+u) = - \Sum_{n=1}^\infty \frac {(-u)^n} {n} = u - \tfrac 1 2 u^2 + \tfrac 1 3 u^3 - \tfrac 1 4 u^4 + \tfrac 1 5 u^5 - \ldots
\end {equation}

if we calculate over $\C$ (or another field of characteristic $\infty$). The corresponding version in characteristic $p$ is obtained as

\begin {equation}
\log_p (1+u) = - \Sum_{n=1}^{p-1} \frac {(-u)^n} {n} = u - \tfrac 1 2 u^2 - \ldots \tfrac 1 {p-2} (-u)^{p-2} - \tfrac 1 {p-1} (-u)^{p-1} \quad \text {(if $u^n = 0$)}
\end {equation}

When dealing with the Lyons group, we have to set $p = 5$; for the convenience of the reader, we communicate the explicit expressions for this particular case, namely

\begin {equation}
\exp_p u = 1 + u + \tfrac 1 2 u^2 + \tfrac 1 6 u^3 + \tfrac 1 {24} u^4 = 1 + u + 3 u^2 + u^3 + 4 u^4
\end {equation}

and

\begin {equation}
\log_p (1+u) = u - \tfrac 1 2 u^2 + \tfrac 1 3 u^3 - \tfrac 1 4 u^4 = u + 2 u^2 + 2 u^3 + u^4
\end {equation}

both being valid only for matrices $u$ over fields of characteristic 5 with $u^5 = 0$.

After these preliminaries we are now ready to plunge into the heart of the matter. We begin with the construction of five matrices in $\F_5^{(111,111)}$ which are the fundamental building bricks for all further calculations.

\begin {definition}

\label {d.mat}

\ \\The symbols $\alpha$, $\beta$, $\gamma$, $\eta$ and $f$ will be reserved for the following matrices:

\begin {enumerate}

\item $\alpha$ is obtained from $[5,8][6,15][7,13][9,12][10,16][11,14]$ by replacing each of the 14 diagonal entries in the positions 6, 7, 8, 9, 11, 12, 13, 17, 20, 21, 24, 25, 26, 27 by $-1 = 4$;

\item $\beta$ is the product of $(4,6,8)(5,7,9)$ and $[2,3,4][5,12,16][6,10,11][7,14,9][8,15,13]$ (in any order);

\item $\gamma$ is calculated with the auxiliary matrices $A = \left| \begin {smallmatrix} 1 \end {smallmatrix} \right|$ and $B = \left| \begin {smallmatrix} . & 1 \\ 4 & 4 \end {smallmatrix} \right|$ via

\begin {equation}
\gamma = [1_3 \otimes A] \oplus [1_3 \otimes B] \oplus [1_3 \otimes (A \oplus A \oplus B \oplus B)] \oplus [1_{12} \otimes (A \oplus B \oplus B \oplus B)]
\end {equation}

\item $\eta$ has exactly 16 nonzero blocks, namely

\begin {eqnarray*}
\eta_{[1,16]} = \begin {array} {|c|} \begin {smallmatrix}
. & . & . & 4 & 3 & . & . \\ . & . & . & . & . & 2 & . \\ 3 & 1 & . & . & . & . & . \\ 2 & 4 & 4 & . & . & . & . \\ 3 & . & 4 & . & . & . & . \\ 1 & 1 & 4 & . & . & 1 & 1 \\ 1 & 2 & 1 & . & . & . & 4 \\ 4 & 4 & 1 & . & . & 1 & 1 \\ 4 & 3 & 4 & . & . & . & 4
\end {smallmatrix} \end {array} \hspace {3mm}
\eta_{[2,13]} = \begin {array} {|c|} \begin {smallmatrix}
2 & 1 & . & 3 & 1 & 4 & . \\ 4 & 2 & . & 4 & 3 & 3 & . \\ . & 3 & 1 & 2 & . & 2 & 4 \\ 4 & 2 & 2 & 2 & 1 & . & 3 \\ . & 2 & 4 & 2 & . & 3 & 1 \\ 1 & 3 & 3 & 2 & 1 & . & 2
\end {smallmatrix} \end {array} \hspace {3mm}
\eta_{[3,7]} = \begin {array} {|c|} \begin {smallmatrix}
3 & 4 & . & 3 & 1 & 4 & . \\ 4 & 2 & . & 1 & 2 & 2 & . \\ 3 & 4 & 4 & 4 & 2 & . & 4 \\ 2 & . & 4 & 2 & 3 & 4 & 4 \\ 3 & . & 1 & 2 & 3 & 1 & 1 \\ . & 4 & 3 & 4 & . & 4 & 3
\end {smallmatrix} \end {array} \hspace {3mm}
\eta_{[4,10]} = \begin {array} {|c|} \begin {smallmatrix}
. & . & . & 2 & 4 & 3 & . \\ . & . & . & 3 & 1 & 3 & . \\ 4 & 3 & 4 & 3 & 4 & . & . \\ 2 & 3 & . & 1 & 2 & . & . \\ 2 & 2 & 3 & 2 & 3 & . & . \\ 2 & 4 & 2 & 1 & 3 & . & .
\end {smallmatrix} \end {array}
\end {eqnarray*}

\begin {eqnarray*}
\eta_{[5,6]} = \begin {array} {|c|} \begin {smallmatrix}
4 & 2 & . & 2 & 2 & 3 & . \\ 4 & 2 & 4 & . & 2 & 3 & 2 \\ 3 & 1 & 3 & 1 & 1 & 3 & 2 \\ 2 & 4 & . & . & 2 & 4 & 4 \\ 1 & . & 2 & 4 & . & . & . \\ 4 & . & 4 & 2 & 3 & . & 3 \\ 3 & . & . & 1 & 4 & 4 & 1
\end {smallmatrix} \end {array} \hspace {3mm}
\eta_{[6,14]} = \begin {array} {|c|} \begin {smallmatrix}
1 & 2 & 4 & 3 & . & 1 & 4 \\ . & 1 & . & 2 & . & 4 & 4 \\ 1 & 4 & 4 & . & . & 2 & 2 \\ 1 & 1 & 3 & 1 & 2 & 1 & 2 \\ 2 & 2 & 1 & 3 & 3 & 2 & 4 \\ 2 & . & 4 & . & . & 1 & . \\ 3 & . & 4 & 1 & . & 4 & 4
\end {smallmatrix} \end {array} \hspace {3mm}
\eta_{[7,2]} = \begin {array} {|c|} \begin {smallmatrix}
1 & 4 & 3 & 4 & 3 \\ 3 & 1 & 2 & 4 & 2 & 4 \\ 1 & 2 & . & 3 & . & 3 \\ . & . & 4 & . & 1 & . \\ 4 & 2 & 4 & 2 & 1 & 3 \\ 2 & 4 & 2 & 4 & 2 & 4 \\ 4 & 3 & 3 & 3 & 3 & 3
\end {smallmatrix} \end {array} \hspace {3mm}
\eta_{[8,12]} = \begin {array} {|c|} \begin {smallmatrix}
4 & 3 & 1 & 3 & . & 1 & 4 \\ . & 4 & . & 2 & . & 4 & 4 \\ 4 & 1 & 1 & . & . & 2 & 2 \\ 1 & 1 & 3 & 4 & 3 & 4 & 3 \\ 2 & 2 & 1 & 2 & 2 & 3 & 1 \\ 2 & . & 4 & . & . & 4 & . \\ 3 & . & 4 & 4 & . & 1 & 1
\end {smallmatrix} \end {array}
\end {eqnarray*}

\begin {eqnarray*}
\eta_{[9,5]} = \begin {array} {|c|} \begin {smallmatrix} 2 & 1 & 2 & 1 & . & 4 & 2 \\ . & 2 & . & 4 & . & . & 2 \\ 3 & 3 & 3 & 1 & . & . & 2 \\ 2 & 1 & 4 & 3 & 1 & 3 & 1 \\ 4 & 2 & 3 & 2 & . & 1 & 2 \\ 4 & 2 & 2 & . & . & 2 & . \\ . & 1 & 1 & 3 & . & 3 & 3
\end {smallmatrix} \end {array} \hspace {3mm}
\eta_{[10,1]} = \begin {array} {|c|} \begin {smallmatrix}
. & . & 4 & 4 & 1 & 1 & 1 & 4 & 4 \\ . & . & 1 & 3 & 2 & 3 & 3 & 2 & 2 \\ . & . & 2 & 2 & . & 3 & 4 & 2 & 1 \\ . & . & . & . & . & . & . & . & . \\ 3 & . & . & . & . & . & . & . & . \\ . & 4 & . & . & . & 1 & 1 & 1 & 1 \\ . & 3 & . & . & . & . & 4 & . & 4
\end {smallmatrix} \end {array} \hspace {3mm}
\eta_{[11,15]} = \begin {array} {|c|} \begin {smallmatrix}
3 & 4 & 3 & 1 & . & 4 & 2 \\ . & 3 & . & 4 & . & . & 2 \\ 2 & 2 & 2 & 1 & . & . & 2 \\ 2 & 1 & 4 & 2 & 4 & 2 & 4 \\ 4 & 2 & 3 & 3 & . & 4 & 3 \\ 4 & 2 & 2 & . & . & 3 & . \\ . & 1 & 1 & 2 & . & 2 & 2
\end {smallmatrix} \end {array} \hspace {3mm}
\eta_{[12,11]} = \begin {array} {|c|} \begin {smallmatrix}
1 & 3 & . & . & 3 & 3 & . \\ 1 & 2 & 4 & 3 & 3 & 1 & 3 \\ 2 & 1 & 3 & 2 & 3 & 3 & 4 \\ 2 & 4 & . & 2 & 2 & . & 4 \\ 2 & 3 & 2 & 2 & 3 & . & 1 \\ 4 & 1 & 1 & 1 & 3 & 3 & 3 \\ 3 & 4 & 4 & 1 & 3 & . & 1
\end {smallmatrix} \end {array}
\end {eqnarray*}

\begin {eqnarray*}
\eta_{[13,3]} = \begin {array} {|c|} \begin {smallmatrix}
2 & 1 & 3 & 2 & 4 & 2 \\ 2 & 1 & 4 & 1 & 2 & 1 \\ 4 & 2 & 1 & . & 4 & 4 \\ . & . & 2 & 2 & . & 2 \\ 4 & 3 & 1 & 2 & 4 & 1 \\ 2 & 1 & 1 & 4 & 3 & 4 \\ 4 & 2 & . & 1 & 1 & .
\end {smallmatrix} \end {array} \hspace {3mm}
\eta_{[14,9]} = \begin {array} {|c|} \begin {smallmatrix}
4 & 2 & . & . & 3 & 3 & . \\ 4 & 3 & 1 & 3 & 3 & 1 & 3 \\ 3 & 4 & 2 & 2 & 3 & 3 & 4 \\ 2 & 4 & . & 3 & 3 & . & 1 \\ 2 & 3 & 2 & 3 & 2 & . & 4 \\ 4 & 1 & 1 & 4 & 2 & 2 & 2 \\ 3 & 4 & 4 & 4 & 2 & . & 4
\end {smallmatrix} \end {array} \hspace {3mm}
\eta_{[15,8]} = \begin {array} {|c|} \begin {smallmatrix}
1 & 3 & . & 2 & 2 & 3 & . \\ 1 & 3 & 1 & . & 2 & 3 & 2 \\ 2 & 4 & 2 & 1 & 1 & 3 & 2 \\ 2 & 4 & . & . & 3 & 1 & 1 \\ 1 & . & 2 & 1 & . & . & . \\ 4 & . & 4 & 3 & 2 & . & 2 \\ 3 & . & . & 4 & 1 & 1 & 4
\end {smallmatrix} \end {array} \hspace {3mm}
\eta_{[16,4]} = \begin {array} {|c|} \begin {smallmatrix}
. & . & . & 3 & 4 & 4 \\ . & . & . & 4 & 2 & 2 \\ . & . & 2 & 4 & 2 & 1 \\ . & . & 1 & 3 & 1 & 4 \\ 1 & 4 & 2 & 4 & 3 & 4 \\ 4 & 4 & . & . & . & . \\ 3 & 3 & . & . & . & .
\end {smallmatrix} \end {array}
\end {eqnarray*}
\begin {equation}
{\ }
\end {equation}

\item $f$ is defined via the following chain of equations: With $D = \left| \begin {smallmatrix} 3 & 1 \\ 1 & 3 \end {smallmatrix} \right|$ and $J = \left| \begin {smallmatrix} . & 1 \\ 1 & . \end {smallmatrix} \right|$ as well as the auxiliary quantities

\begin {equation}
F_1 = \op {diag} (4,4,3) \oplus D \oplus D \oplus D
\end {equation}

\begin {equation}
F_2 = 1_3 \otimes [1_2 \oplus (2D) \oplus(-2D)]
\end {equation}

\begin {equation}
F_4 = J \otimes 1_6 \otimes [1_1 \oplus (-D) \oplus D \oplus (-D)]
\end {equation}

we put

\begin {equation}
f = F_1 \oplus F_2 \oplus F_4
\end {equation}

\end {enumerate}

\end {definition}

\begin {remarks}

\begin {enumerate}

\item As already noticed before, all (111,111)-dimensional matrices in the remainder of this paper will be expressed in terms of $\alpha$, $\beta$, $\gamma$, $\eta$ and $f$.

\item The five matrices defined here are block-monomial, but for only one of them ($\eta$) is it necessary to give the nonzero blocks explicitly.

\item In fact, $\alpha$, $\beta$, $\gamma$ and $\eta$ are required (and suffice) to construct the group $X$ which realizes the absolute minimal representation of the sporadic Lyons group $Ly$, while $f$ is merely needed for the construction of an $X$-invariant quadratic form.

\item Since $A$ is the 1-dimensional identity matrix, the defining formula for $\gamma$ can be simplified to

\begin {equation}
\gamma = 1_3 \oplus [1_3 \otimes B] \oplus [1_3 \otimes (1_2 \oplus B \oplus B)] \oplus [1_{12} \otimes (1_1 \oplus B \oplus B \oplus B)]
\end {equation}

\end {enumerate}

\end {remarks}

The relations $\alpha^2 = \beta^3 =1$, $\alpha^\beta = \omega$ and $\omega^\beta = \alpha \omega = \omega \alpha$ form a presentation of the alternating group on 4 letters. Since $\alpha \ne 1$ we more precisely have $\la\alpha,\beta\ra \iso A_4$. Furthermore, $\gamma$ is of order 3 and commutes with $\alpha$ and $\beta$. Hence

\begin {equation}
K = \la\alpha, \beta, \gamma\ra \iso \ol K = A_{\{1,2,3,4\}} \times A_{\{a,b,c\}} \iso A_4 \times A_3
\end {equation}

An explicit isomorphism from $K$ onto $\ol K$ which we shall fix from now on is given by

\begin {eqnarray}
\ol \alpha &= (1,2)(3,4)\\
\ol \beta &= (1,2,3) \\
\ol \gamma &= (a,b,c)
\end {eqnarray}

We use this to transfer the natural action of $\ol K$ on (the names of) the points, lines, planes and flags of Kantor's apartment to $K$. For instance, the element $\alpha \beta \gamma \in K$ corresponds to

\begin {equation}
\ol {\alpha \beta \gamma} = (1,2)(3,4) \cdot (1,2,3) \cdot (a,b,c) = (1,3,4)(a,b,c)
\end {equation}

and thus maps the line $(1a,2b)$ to

\begin {equation}
(1a,2b)^{\alpha \beta \gamma} = (1a,2b)^{(1,3,4)(a,b,c)} = (3b,2c)
\end {equation}

\begin {remark}

This obviously provides us with an action of $K$ (or $\ol K$) on the apartment which is transitive on the points as well as sharply transitive and orientation-preserving on the lines.

\end {remark}

The block permutation part of $\eta$, namely

\begin {equation}
(1,16,4,10)(2,13,3,7)(5,6,14,9)(8,12,11,15)
\end {equation}

has order 4, hence $\eta^4$ is block-diagonal. In fact, 4 of the 16 diagonal blocks of $\eta^4$ vanish (those in positions 2, 3, 10 and 16) while the other 12 are non-zero. Thus, $\eta^4 \ne 0$. Multiplication of $\eta^4$ with $\eta$, however, shows that $\eta^5$ is the zero matrix. Hence $\eta$ is nilpotent of order 5:

\begin {equation}
\eta^4 \ne \eta^5 = 0
\end {equation}

and we therefore may apply the exponential function to $\eta$. We set

\begin {equation}
\xi = \exp_5(\eta) = \Sum_{m=0}^4 \frac {\eta^m} {m!} = 1 + \eta + \tfrac 1 2 \eta^2 + \tfrac 1 6 \eta^3 + \tfrac 1 {24} \eta^4 = 1 + \eta + 3 \eta^2 + \eta^3 + 4 \eta^4
\end {equation}

which is of the general shape

\begin {equation}
\xi - 1 = \eta + \rm O(\eta^2)
\end {equation}

By induction we easily get

\begin {equation}
(\xi-1)^n = \eta^n + \rm O(\eta^{n+1})
\end {equation}

for all $n \in \N$ and in particular with $n = 4$ and $n =5$:

\begin {eqnarray}
(\xi-1)^4 = \eta^4 \ne 0 \\
(\xi-1)^5 = \eta^5 = 0
\end {eqnarray}

The penultimate result implies $\xi \ne 1$, while the last equation can be simplified considerably since we calculate modulo 5. This leads to

\begin {equation}
0 = (\xi-1)^5 = \xi^5 - 5 \xi^4 + 10 \xi^3 - 10 \xi^2 + 5 \xi - 1 = \xi^5 - 1
\end {equation}

Consequently, $\xi$ has order 5 and is thus invertible. The Frobenius automorphism of the basic field $\F_5$ provides us with a slightly more definite variant, namely

\begin {equation}
\det \xi = (\det \xi)^5 = \det \xi^5 = \det 1 =1
\end {equation}

which is tantamount to $\xi \in SL_{111}(5)$.

We now employ the fact that the action of $K$ on $\frak L$ is simply transitive and orientation-preserving to construct further matrices and matrix groups. \\

\begin {definition}

\label {d.x}

\ \\Let $L$ be any line in $\frak L$. We determine the unique $k \in K$ with $L = (1a,2b)^k = (1a,2b)^{\ol k}$ and thereafter set

\begin {equation}
x_L = \xi^k = k^{-1} \xi k \in SL_{111}(5)
\end {equation}

\begin {equation}
X_L = \grp {x_L} \iso 5
\end {equation}

and finally

\begin {equation}
X = \gruppe {X_L} {L \in \frak L} = \gruppe {x_L} {L \in \frak L} \le SL_{111}(5)
\end {equation}

\end {definition}

\section {Configurations and configuration groups}

\label {a.conf}

The main purpose of the present paper is to study the just defined group $X$. To that end, we first have to investigate certain natural geometric subgroups which we want to introduce next.

\begin {definition}

\label {d.conf}

\ \\Consider an oriented line $L = (P,Q) \in \frak L$.

\begin {enumerate}

\item We associate with $L$ the following {\bf configurations} (\,= subsets of $\frak L$):

\begin {enumerate}

\item The {\bf star} $\mathcal S(L)$, consisting of the six lines $\menge {L_i} {i \in \F_7^\times}$ incident with $P$. We set $L_1 =L$ and name the other lines such that $L_i$ under the $60^\circ$-rotation with center $P$ is mapped to $- L_{3i}$, where the minus sign denotes orientation reversal. The unique point other than $P$ inciding with $L_i$ is called $Q_i$ (Fig. \ref {abb.s}).

\item The {\bf hexagon} $\mathcal H(L)$, composed of the $L_i$ and six further lines $l_i$, $i \in \F_7^\times$, where $l_i$ incides with the points $Q_{2i}$ and $Q_{3i}$ (Fig. \ref {abb.h}).

\item The {\bf quartet} $\mathcal Q(L)$, containing the lines $l_1$ and $l_6$ together with their direct continuations $m_1$ resp. $m_6$ (Fig. \ref {abb.q}).

\end {enumerate}

\item For the sake of brevity we usually replace the somewhat cumbersome names $x_{L_i}$, $x_{l_i}$ and $x_{m_i}$ for the $X$-generators associated with the lines $L_i$, $l_i$ resp. $m_i$ by $\Lambda_i(L)$, $\lambda_i(L)$ and $\mu_i(L)$. If it is evident from the context to which line $L$ we refer, we may omit the argument and simply write $\Lambda_i$, $\lambda_i$ or $\mu_i$.

\item Furthermore we define the {\bf configuration groups} with respect to the line $L$ (which is always displayed in red):

\end {enumerate}

\begin {enumerate}

\item {\bf Star group}

\begin {equation}
\frak S(L) = \gruppe {\Lambda_i} {i \in \F_7^\times}
\end {equation}

\item {\bf Hexagon group}

\begin {equation}
\frak H(L) = \gruppe {\Lambda_i, \lambda_i} {i \in \F_7^\times}
\end {equation}

\item {\bf Quartet group}

\begin {equation}
\frak Q(L) = \grp {\lambda_1, \mu_1, \lambda_6, \mu_6}
\end {equation}

\end {enumerate}

\end {definition}

\begin {remark}

The line $L$ itself is (as $L_1$) part of $\mathcal S(L)$ and $\mathcal H(L)$, but not of $\mathcal Q(L)$. This is the reason why in the graph corresponding to the latter case $L$ is represented by a dashed line.

\end {remark}

\begin {center}

\abb {5.1cm} {5.1cm} {star} {\bf Star $\mathcal S(L)$)} \label {abb.s} \abb {5.1cm} {5.1cm} {hexagon} {\bf Hexagon $\mathcal H(L)$} \label {abb.h} \abb {5.1cm} {5.1cm} {quartet} {\bf Quartet $\mathcal Q(L)$} \label {abb.q}

\end {center}

We now come to an obvious but very useful principle:

\begin {lemma}

\label {l.wort}

\ \\ For $L, L' \in \frak L$ let $w$ be some word in the generators $\Lambda_i(L), \lambda_j(L), \mu_l(L)$ and $w'$ the identically built word in $\Lambda_i(L'), \lambda_j(L'), \mu_l(L')$. Then $w$ and $w'$ are conjugate in $X$. In particular, the equations $w = 1$ and $w' = 1$ are tantamount.

\end {lemma}

\begin {proof}

By construction, there is a unique $k \in K \le X$ with $L^k = L'$ (see the above remark). Conjugation by $k$ thus transforms $w$ into $w'$, and the result immediately follows.

\end {proof}

Our next aim is to prove several properties of the configuration groups. \\

\begin {lemma}

\label {l.hexrel}

\ \\ Assume $L \in \frak L$. Let $i$ be a square in $\F_7^\times$; in other words $i \in \{1,2,4\}$. Then the following relations hold in the hexagon group $\frak H(L)$ associated with $L$:

\begin {align}
\left[ \Lambda_i,\Lambda_{-4i} \right] &= 1 \\
\left[ \Lambda_i,\Lambda_{2i} \right] &= \Lambda_{-4i}^4 \\
\left[ \lambda_i,\lambda_{-4i} \right] &= \Lambda_2 \\
\left[ \lambda_i,\lambda_2 \right] &= \Lambda_2 \lambda_{-4i}^3 \Lambda_{-i}^3 \\
\left[ \Lambda_i,\lambda_i \right] &= 1 \\
\left[ \Lambda_i,\lambda_{-4i} \right] &= \lambda_{-2i}^3 \Lambda_{-4i} \lambda_i^2 \Lambda_2^4 \\
\left[ \Lambda_i,\lambda_2 \right] &= \lambda_{4i}^2 \Lambda_{-2i} \lambda_{-i}^2 \Lambda_{4i} \\
\left[ \Lambda_i,\lambda_{-i} \right] &= 1 \\
\left[ \Lambda_i,\lambda_{4i} \right] &= 1 \\
\left[ \Lambda_i,\lambda_{-2i} \right] &= 1
\end {align}

\begin {align}
\left[ \Lambda_{-i},\Lambda_{4i} \right] &= 1 \\
\left[ \Lambda_{-i},\Lambda_{-2i} \right] &= \Lambda_{4i}^3 \\
\left[ \lambda_{-i},\lambda_{4i} \right] &= \Lambda_{-2i}^2 \\
\left[ \lambda_{-i},\lambda_{-2i} \right] &= \Lambda_{-2i}^2 \lambda_{4i}^3 \Lambda_i^4 \\
\left[ \Lambda_{-i},\lambda_{-i} \right] &= 1 \\
\left[ \Lambda_{-i},\lambda_{4i} \right] &= \lambda_2^4 \Lambda_{4i}^2 \lambda_{-i} \Lambda_{-2i}^4 \\
\left[ \Lambda_{-i},\lambda_{-2i} \right] &= \lambda_{-4i} \Lambda_2^2 \lambda_i \Lambda_{-4i} \\
\left[ \Lambda_{-i},\lambda_i \right] &= 1 \\
\left[ \Lambda_{-i},\lambda_{-4i} \right] &= 1 \\
\left[ \Lambda_{-i},\lambda_2 \right] &= 1
\end {align}

Furthermore, for both groups $\la \Lambda_i, \Lambda_{-i} \ra$ and $\la \lambda_i, \lambda_{-i} \ra$, there exists an isomorphism to $SL_2(5)$ such that one generator is mapped to an upper and the other to a lower triangular matrix.

\end {lemma}

\begin {proof}

By Lemma \ref {l.wort} it suffices to prove the assertions for the special case $L = (1a, 2b)$. Moreover, we only have to consider $i = 1$, since all other relations then follow by conjugation with powers of $(2,3,4) \in K$. This reduces the required amount of labour by a factor of $36 \cdot 3 = 108$. The remaining formulas are then quickly verified, and the first part of the lemma is proved.

The second proposition is also easily established. Explicit isomorphisms are, for example,

\begin {equation}
\Lambda_i \longmapsto \left| \begin {matrix} 1 & 1 \\ . & 1 \end {matrix} \right| \qquad \Lambda_{-i} \longmapsto \left| \begin {matrix} 1 & . \\ 3 & 1 \end {matrix} \right|
\end {equation}

and

\begin {equation}
\lambda_i \longmapsto \left| \begin {matrix} 1 & 1 \\ . & 1 \end {matrix} \right| \qquad \lambda_{-i} \longmapsto \left| \begin {matrix} 1 & . \\ 4 & 1 \end {matrix} \right|
\end {equation}

in the respective cases.

\end {proof}

This result enables us to characterize both the hexagon and the star groups.

\begin {lemma}

\label {l.hex}

Let $L \in \frak L$. Then

\begin {enumerate}

\item $\frak H(L) \iso G_2(5)$;

\item $\frak S(L) \iso A_2(5)$;

\item in particular, $\frak H(L)$ and $\frak S(L)$ are simple Chevalley groups of orders

\begin {equation}
|\frak H(L)| = 2^6 \cdot 3^3 \cdot 5^6 \cdot 7 \cdot 31 = 5859000000
\end {equation}

\begin {equation}
|\frak S(L)| = 2^5 \cdot 3 \cdot 5^3 \cdot 31 = 372000
\end {equation}

\item Furthermore, the map $x \longmapsto x^*$ with

\begin {equation}
\Lambda^*_1 = \left| \begin {matrix} 1 & . & . \\ . & 1 & 3 \\ . & . & 1 \end {matrix} \right| \qquad \quad
\Lambda^*_2 = \left| \begin {matrix} 1 & . & . \\ . & 1 & . \\ 3 & . & 1 \end {matrix} \right| \qquad \quad
\Lambda^*_4 = \left| \begin {matrix} 1 & 3 & . \\ . & 1 & . \\ . & . & 1 \end {matrix} \right|
\end {equation}

\begin {equation}
\Lambda^*_6 = \left| \begin {matrix} 1 & . & . \\ . & 1 & . \\ . & 1 & 1 \end {matrix} \right| \qquad \quad
\Lambda^*_5 = \left| \begin {matrix} 1 & . & 1 \\ . & 1 & . \\ . & . & 1 \end {matrix} \right| \qquad \quad
\Lambda^*_3 = \left| \begin {matrix} 1 & . & . \\ 1 & 1 & . \\ . & . & 1 \end {matrix} \right|
\end {equation}

can be extended in a unique way to an isomorphism from $\frak S(L) \iso A_2(5)$ to $SL_3(5)$. We call it the {\bf star isomorphism} (associated with $L$) and keep this notation for the remainder of the paper.

\end {enumerate}

\end {lemma}

\begin {proof}

According to \zitSt, cf. also \zitH or \zitC, the relations given in Lemma \ref {l.hexrel} form a presentation of $G_2(5)$; hence $\frak H(L)$ is a (nontrivial) homomorphic image of that (simple) group and (1) follows.

Moreover, restricting to the formulas containing only the generators of $\frak S(L)$, i.\,e. omitting those referring to the $\lambda_i$, we are left with a presentation of $A_2(5)$, and an analogous argument as before provides us with (2). Proposition (3) then is immediate from Chevalley theory.

Finally we check that the images $\Lambda^*_i$ satisfy all defining relations for $\frak S(L)$ as an abstract group with the generators $\Lambda_i$, and applying the above trick once more leads to assertion (4).

\end {proof}

Our next goal is to determine the structure of the quartet groups.

\begin {lemma}

\label {l.quad}

All quartet groups $\frak Q(L)$, $L \in \frak L$, are of the form

\begin {equation}
\frak Q(L) \iso 2 \Schur A_6 \iso SL_2(9)
\end {equation}

\end {lemma}

\begin {proof}

Thanks to Lemma \ref {l.wort} we only need to consider the special case $L = (1a,2b)$. Omitting the argument $L$, we see that $\mathcal S(L) = \{l_1, m_1, l_6, m_6\}$. As Fig. \ref {abb.q} reveals, the first two of them are contained in $\mathcal S(m_1)$:

\begin {equation}
l_1 = l_1(L) = L_6(m_1) = (4b,3c) \qquad m_1 = m_1(L) = L_1(m_1) = (3c,4a)
\end {equation}

while the other two can be likewise interpreted as elements of $\mathcal S(m_6)$:

\begin {equation}
l_6 = l_6(L) = L_6(m_6) = (3b,4c) \qquad m_6 = m_6(L) = L_1(m_6) = (4c,3a)
\end {equation}

Let $s \in \{1,6\}$. The images of $\lambda_s$ and $\mu_s$ under the star isomorphism from $\mathcal S(m_s)$ onto $SL_3(5)$ introduced in Lemma \ref {l.hex}, namely

\begin {equation}
\lambda^*_s =
\left| \begin {matrix}
1 & . & . \\
. & 1 & . \\
. & 1 & 1
\end {matrix} \right|
\qquad
\mu^*_s =
\left| \begin {matrix}
1 & . & . \\
. & 1 & 3 \\
. & . & 1
\end {matrix} \right|
\end {equation}

are independent of $s$ and obviously generate a subgroup isomorphic to $SL_2(5)$. We define the auxiliary elements

\begin {align}
h_1 &= \mu_1 \lambda_1^4 \mu_1^4 \lambda_1^3 \\
h_2 &= \mu_1^4 \lambda_1 \\
h_3 &= \lambda_6^2 \mu_6^4 \lambda_6^4 \\
h_4 &= \lambda_6^2 \mu_6^2
\end {align}

and set $H = \la h_1, h_2, h_3, h_4 \ra$. Since, as may be verified easily, these relations can be solved for $\lambda_s$ and $\mu_s$, e.\,g. via

\begin {align}
\lambda_1 &= h_1 h_4^2 h_2^2 \\
\mu_1 &= h_1 h_4^2 h_2 \\
\lambda_6 &= h_2^2 h_4 h_3^2 \\
\mu_6 &= h_3 h_2 h_4 h_1 h_3 h_2^2 h_4 h_1
\end {align}

we obtain

\begin {equation}
\frak Q(L) = \la \lambda_1, \mu_1, \lambda_6, \mu_6 \ra = \la h_1, h_2, h_3, h_4 \ra = H
\end {equation}

Hence it merely remains to prove $H \iso 2 \Schur A_6 \iso SL_2(9)$.

For all $i,j \in \{1,2,3,4\}$ with $i \ne j$ the expressions $(h_i h_j)^2$ have the same value, while $h_i^3 = 1$. This is the famous presentation for the Schur double cover of $A_6$ discovered by \zitMo.

Consequently, $H$ is isomorphic to some factor group of $2 \Schur A_6$, namely $2 \Schur A_6$, $A_6$ or $1$. Since $\frak Q(L)$ possesses $SL_2(5)$-subgroups which contain elements of order 10, while neither $A_6$ nor $1$ does, we conclude $H \iso 2 \Schur A_6$.

An explicit isomorphism of $H$ with $SL_2(9)$ is defined by

\begin {equation}
\{h_1,h_2,h_3,h_4\} \longmapsto \mng {\left| \begin {matrix} . & 1 \\ -1 & -1 \end {matrix} \right|, \quad
\left| \begin {matrix} . & -1 \\ 1 & -1 \end {matrix} \right|, \quad
\left| \begin {matrix} -1 & i \\ i & . \end {matrix} \right|, \quad
\left| \begin {matrix} -1 & -i \\ -i & . \end {matrix} \right|}
\end {equation}

where the matrix entries are elements of $\F_9$ and $i$ denotes any of the two solutions of $i^2 + 1 = 0$ in $\F_9$. The proof is entirely analogous to the above argument.

\end {proof}

\section {Identification of $X$; apartment symmetries}

\label {a.idx}

\ \\ We are now prepared to determine the structure of $X$.

\begin {theorem}

\label {s.ly}

{\ }

\begin {enumerate}

\item $X$ is isomorphic with the sporadic group $Ly$ of Lyons.

\item In particular, $X$ is simple of order

\begin {equation}
|Ly| = 2^8 \cdot 3^7 \cdot 5^6 \cdot 7 \cdot 11 \cdot 31 \cdot 37 \cdot 67 = 51765179004000000
\end {equation}

\item $X$ realizes an absolute minimal representation of the Lyons group (111-dimensional over $\F_5$).

\end {enumerate}

\end {theorem}

\begin {proof}

By the main result (Satz 3.4.2) of \zitG, Lemma \ref {l.hex} (first part) and Lemma \ref {l.quad} together imply (1); assertion (2) then follows from the characterization of $Ly$ used in the same paper.

Finally, statement (3) is an immediate consequence of the modular character theoretic investigations in \zitMN.

\end {proof}

Next we want to investigate the geometry of the Kantor apartment in some detail.

\begin {definition}

\label {d.w}

\ \\ The {\bf Weyl group} of $X$ (or $\frak A$) is

\begin {equation}
W = S_{\{1,2,3,4\}} \times S_{\{{a,b,c}\}} \iso S_4 \times S_3
\end {equation}

acting in the natural manner on the (names of the) points, lines, planes and flags of $\frak A$.

\end {definition}

Even at first glance $\frak A$ exhibits numerous symmetries. Let us illustrate this with a few typical examples:

\begin{remark}

\begin {enumerate}

\item Given two (not necessarily distinct) points $P,Q \in \frak P$, there is a unique {\bf translation} from $P$ to $Q$, namely the Weyl element which shifts $\frak A$ parallel to itself and maps $P$ to $Q$. For instance, if $P = 1a$ and $Q = 4b$, we get $(1,4)(2,3)(a,b,c) \in W$. All 12 translations together form the abelian subgroup $[A_4,A_4] \times A_3 \iso 2^2 \times 3$ of $W$.

\item The left-right reflection of the diagram in Fig. \ref {abb.k} interchanges the points incident with the line $(3b,4c)$ or equivalently those inciding with $(4b,3c)$. This amounts to $(3,4)(b,c) \in W$.

\item The up-down reflection as shown in Fig. \ref {abb.k} is represented by $(3,4) \in W$.

\item The (counterclockwise) $60^\circ$-rotation around the central point $1a$ in the same diagram is $(2,3,4)(b,c) \in W$.

\end {enumerate}

Note that the transformations under (1) and (3) are orientation-preserving, while those under (2) and (4) are orientation-reversing. In particular, the two types of mirror symmetries are inequivalent to each other.

\end {remark}

All automorphisms of $\frak A$ (as a simplicial complex) are now described easily.

\begin {lemma}

\label {l.w}

\ \\ The full automorphism group $\aut \frak A$ of Kantor's (unoriented) apartment $\frak A$ is identical with the Weyl group. A Weyl element preserves or reverses the orientation of $\frak A$ (cf. Fig. \ref {abb.nm}) according as it is contained in the naturally embedded subgroup $S_4 \times A_3$ of $W$ or not. Moreover, $W$ acts sharply transitively on the 144 maximal flags in $\frak A$.

\end {lemma}

\begin {proof}

Obvious.

\end {proof}

It is now possible to introduce some further interesting elements and a very important subgroup of $X$.

\begin {definition}

\label {d.r}

\ \\ For each $L \in \frak L$ we omit, as usual, the argument $L$ of $\Lambda_i$ and $\lambda_i$ and set

\begin {align}
R(L) &= \Lambda_1^2 \Lambda_6^4 \Lambda_1^2 \\
r(L) &= \lambda_1^2 \lambda_6^3 \lambda_1^2
\end {align}

Under the isomorphisms from $\la \Lambda_1, \Lambda_6 \ra$ and $\la \lambda_1, \lambda_6 \ra$ to $SL_2(5)$ described in the proof of Lemma \ref {l.hexrel}, the elements $R(L)$ and $r(L)$ each correspond to the matrix $\left| \begin {smallmatrix} . & 2 \\ 2 & . \end {smallmatrix} \right|$ and act (by conjugation) on $\frak H(L)$ as horizontal and vertical reflections. Moreover we define the group

\begin {equation}
N = \gruppe {R(L), r(L)} {L \in \frak L}
\end {equation}

\end {definition}

For most applications, the large number of generators in the definitions of $N$ and $X$ is quite inconvenient. Therefore it is advisable to replace them by smaller sets. As concerns $N$, a suitable choice is given by

\begin {lemma}

\label {l.Nerz}

\ \\ $N$ is generated by the four elements

\begin {align}
n_1 &= \beta \\
n_2 &= x_{(3a,4b)} x_{(4a,3b)} x_{(3a,4b)} \\
n_3 &= \gamma \\
n_4 &= x_{(1b,2c)}^3 x_{(2c,1a)} x_{(1b,2c)}^3 x_{(3a,4b)} x_{(4a,3b)} x_{(3a,4b)}
\end {align}

whose images in the Weyl group $W \iso \ol N = N/T$ are

\begin {align}
\ol {\rule{0pt}{8pt} n_1} &= (1,2,3) \\
\ol {\rule{0pt}{8pt} n_2} &= (3,4) \\
\ol {\rule{0pt}{8pt} n_3} &= (a,b,c) \\
\ol {\rule{0pt}{8pt} n_4} &= (a,b)
\end {align}

\end {lemma}

\begin {proof}

The $W$-images of $n_1, \ldots, n_4$ obviously generate the whole Weyl group; hence it suffices to show $T \le \la n_1,n_2,n_3,n_4 \ra$. In fact, we find by trial and error that even

\begin {equation}
T = \la n_1 n_4 n_1^{-1} n_4, n_1^{-1} n_4 n_1 n_4 \ra \le \la n_1,n_4\ra
\end {equation}

holds.

\end {proof}

In a similar vein we have

\begin {lemma}

\label {l.erz}

\ \\ $X = \la \alpha, \beta, \gamma, \xi \ra$.

\end {lemma}

\begin {proof}

$Y = \la \alpha, \beta, \gamma, \xi \ra$ contains $\la \alpha, \beta, \gamma \ra = K$, consequently also $\xi^K = \menge {\xi^k} {k \in K} = \menge {x_L} {L \in \frak L}$ and the group $X$ spanned by these elements.

For a proof of the reverse inclusion $Y \le X$ we need only show that the generators of $Y$ lie in $X$. For $\xi = x_{(1a,2b)}$ this is obvious; for the other three matrices it follows e.\,g. from the easily verified equations

\begin {align}
\alpha &= R_{(3a,1b)} R_{(1a,3b)} R_{(3a,4b)} R_{(2a,1b)}^{-1} \\
\beta &= R_{(1a,2b)} R_{(3a,4b)}^{-1} R_{(2a,1b)} R_{(4a,2b)}^{-1} \\
\gamma &= R_{(1b,2c)} R_{(1a,2b)}^{-1}
\end {align}

\end {proof}

We give an application of the last lemma.

\begin {theorem}

\label {s.sp}

\ \\ The symmetric {\bf scalar product} $\la \pkp \ra: V \times V \longrightarrow \F_5$ given by

\begin {equation}
\la u,v \ra = u \cdot f \cdot \tr v
\end {equation}

on $V$ and the associated {\bf quadratic form} $\mathcal F: V \longrightarrow \F_5$ with

\begin {equation}
\mathcal F(u) = \la u,u \ra
\end {equation}

are $X$-invariant and non-degenerate.

\end {theorem}

\begin {proof}

Both invariance statements are equivalent to the relation $x f \, \tr x = f$ for all $x \in X$. Clearly, this condition has only to be verified for $x$ in a generating subset of $X$, for instance $\{ \alpha, \beta, \gamma, \xi \}$, which is easily done. The nondegeneracy is tantamount to the non-vanishing of the determinant of $f$ which in fact has the value $4 = -1$, as can be derived directly from the definition.

\end {proof}

\section {Kantor's standard torus}

\label {a.tor}

Kantor's theory of the Lyons group and particularly its 5-geometry is essentially based on a certain torus (\zitK). The original approach was inspired by the properties of the homonymous objects in Chevalley groups, and therefore tori were defined as complements of Sylow 5-subgroups in their normalizers (the Borel groups). For our present purpose, however, a simpler characterization is availabe, since due to \zitNM for the relevant groups $A_2(5)$, $G_2(5)$ and $Ly$ itself, a torus is just an arbitrary subgroup isomorphic to $4^2$.

This makes it very easy to construct tori galore in any of the just mentioned overgroups, but since we also want to establish the associated set of root subgroups we have to be a bit more careful and construct a special torus in a deliberate manner. To achieve that more ambitious goal we first consider the situation in a star group where a natural choice suggests itself:

\begin {definition}

\label {d.lt}

\ \\ To each $L \in \frak L$ we associate the {\bf line torus} $T(L) \le \frak S(L)$ as the group which corresponds via the star isomorphism $x \longmapsto x^*$ to the diagonal subgroup of $SL_3(5)$:

\begin {equation}
T^*(L) = [T(L)]^* = \rm {Diag}(SL_3(5)) = \menge {diag(\rho, \sigma, \tau)} {\rho, \sigma, \tau \in \F_5^\times, \rho \sigma \tau = 1} \iso 4^2
\end {equation}

\end {definition}

We collect some more detailed and useful information about the line tori.

\begin {lemma}

\label {l.lt}

\ \\ Consider $L \in \frak L$ and the corresponding star group $\frak S(L)$. Let $i \in \F_7^\times$. Then

\begin {enumerate}

\item All elements

\begin {equation}
t_i = t_i(L) = \Lambda_i \Lambda_{-i}^2 \Lambda_i^2 \Lambda_{-i}
\end {equation}

have order 4 and are contained in $T(L)$.

\item More explicitly, under the star isomorphism they are mapped as follows:

\begin {eqnarray}
t^*_1 = {\rm diag}(1,2,3) \qquad t^*_2 = {\rm diag}(3,1,2) \qquad t^*_4 = {\rm diag}(2,3,1) \\
t^*_6 = {\rm diag}(1,3,2) \qquad t^*_5 = {\rm diag}(2,1,3) \qquad t^*_3 = {\rm diag}(3,2,1)
\end {eqnarray}

\item Among the identities obeyed by the $t_i$ we note

\begin {equation}
t_{-i} = t_i^{-1}
\end {equation}

while the product of the $t_i$ where $i$ runs over all squares in $\F_7^\times$, namely $i \in \{ 1,2,4 \}$, evaluates to the unit matrix:

\begin {equation}
t_1 t_2 t_4 = 1
\end {equation}

\item $T(L) \iso 4^2$ is generated by any pair $t_i(L)$ and $t_j(L)$, $i,j \in \F_7^\times$, provided $j \ne \pm i$.

\end {enumerate}

\end {lemma}

\begin {proof}

After translation via the star isomorphism, all propositions boil down to simple matrix equations in $SL_3(5)$ and can be checked by hand.

\end {proof}

The situation is much nicer than one might expect.

\begin {theorem}

\label {s.tor}

\ \\ All line tori $T(L)$, $L \in \frak L$, are equal. We henceforth write $T$ instead of $T(L)$ and call $T$ the {\bf Kantor (standard) torus} of $X$. The eigenspaces of $T$ corresponding to its 16 irreducible 5-characters are $E_1, \ldots, E_{16}$.

\end {theorem}

\begin {proof}

Set $L = (1a,2b) \in \frak L$ and $t = t_i(L)$ with $i \in \{ 1,2 \}$. Then $t$ is mapped via the star isomorphism for $\frak S(L)$ to some diagonal matrix $t^* \in SL_3(5)$. As a simple evaluation of the defining formula reveals, the 111-dimensional matrix $t$ itself is also diagonal and -- even more so -- block scalar.

We display $t_1(L)$ and $t_2(L)$ in Table \ref {t.t1t2}. Block scalar matrices are completely determined by the 16 multipliers on $E_1, \ldots, E_{16}$, so only these are shown.

\begin {center}

\begin {tabelle}

\label {t.t1t2}

\begin {tabular} {||c||c|c|c|c|c|c|c|c|c|c|c|c|c|c|c|c||} \hline
\rule {0mm} {4.4mm} & $E_1$ & $E_2$ & $E_3$ & $E_4$ & $E_5$ & $E_6$ & $E_7$ & $E_8$ & $E_9$ & $E_{10}$ & $E_{11}$ & $E_{12}$ & $E_{13}$ & $E_{14}$ & $E_{15}$ & $E_{16}$ \\ \hline
$t_1(L)$ \rule {0mm} {4.4mm} & 1 & 1 & 4 & 4 & 1 & 3 & 2 & 4 & 2 & 2 & 1 & 2 & 3 & 4 & 3 & 3 \\ \hline
$t_2(L)$ \rule {0mm} {4.4mm} & 1 & 4 & 1 & 4 & 2 & 1 & 3 & 2 & 4 & 2 & 3 & 1 & 2 & 3 & 4 & 3 \\ \hline
\end {tabular}

\end {tabelle}

\end {center}

All 16 possible combinations of eigenvalues occur; hence every character is represented by $T(L)$, and each $E_m$ is a full eigenspace. Thus the second part of the theorem holds.

Let now $k \in \{ \alpha, \beta, \gamma \}$. Then $k$ is block monomial and has the canonical decomposition $k = k_D k_P$, where the first factor is block diagonal and the other a block permutation. We find

\begin {equation}
t^k = t^{k_D k_P} = t^{k_P}
\end {equation}

as block scalar and block diagonal matrices always commute.

Consequently $t^k$ is block scalar as well. The eigenvalues of $t^k$ are obtained from those of $t$ by the block permutation $k_P$. They are collected in Table \ref {t.thochk}, continuing the notation of Table \ref {t.t1t2}. For the sake of brevity and readability, here the common argument of $t_1(L)$ and $t_2(L)$ is omitted.

In addition we give in the last column a representation of $t^k$ as a word in the generators of $T(L)$. In all 6 cases $t^k \in T(L)$.

\begin {center}

\begin {tabelle}

\label {t.thochk}

\begin {tabular} {||c||c|c|c|c|c|c|c|c|c|c|c|c|c|c|c|c||c||} \hline
$t^k$ \rule {0mm} {4.4mm} & $E_1$ & $E_2$ & $E_3$ & $E_4$ & $E_5$ & $E_6$ & $E_7$ & $E_8$ & $E_9$ & $E_{10}$ & $E_{11}$ & $E_{12}$ & $E_{13}$ & $E_{14}$ & $E_{15}$ & $E_{16}$ & $t^k$ \\ \hline
$t_1^\alpha$ \rule {0mm} {4.4mm} & 1 & 1 & 4 & 4 & 4 & 3 & 3 & 1 & 2 & 3 & 4 & 2 & 2 & 1 & 3 & 2 & $t_1 t_2^2$ \\ \hline
$t_2^\alpha$ \rule {0mm} {4.4mm} & 1 & 4 & 1 & 4 & 2 & 4 & 2 & 2 & 1 & 3 & 3 & 4 & 3 & 3 & 1 & 2 & $t_1^2 t_2$ \\ \hline
$t_1^\beta$ \rule {0mm} {4.4mm} & 1 & 4 & 1 & 4 & 3 & 1 & 2 & 3 & 4 & 3 & 2 & 1 & 3 & 2 & 4 & 2 & $t_2^3$ \\ \hline
$t_2^\beta$ \rule {0mm} {4.4mm} & 1 & 4 & 4 & 1 & 3 & 3 & 4 & 2 & 3 & 1 & 2 & 2 & 4 & 3 & 2 & 1 & $t_1 t_2^3$ \\ \hline
$t_1^\gamma$ \rule {0mm} {4.4mm} & 1 & 1 & 4 & 4 & 1 & 3 & 2 & 4 & 2 & 2 & 1 & 2 & 3 & 4 & 3 & 3 & $t_1$ \\ \hline
$t_2^\gamma$ \rule {0mm} {4.4mm} & 1 & 4 & 1 & 4 & 2 & 1 & 3 & 2 & 4 & 2 & 3 & 1 & 2 & 3 & 4 & 3 & $t_2$ \\ \hline
\end {tabular}

\end {tabelle}

\end {center}

Applying Lemma \ref {l.wort}, this provides us with

\begin {equation}
T(L^k) = [T(L)]^k = T(L)
\end {equation}

for $k$ as given and {\it a fortiori} for $k \in \la \alpha, \beta, \gamma \ra = K$. This proves the first assertion because of the transitivity of $K$ on $\frak L$.

\end {proof}

\section {The Lyons geometry}

\label {a.lygeo}

Our next task is to clarify some local properties of $X \iso Ly$. We start with the centralizer of the standard torus.

\begin {lemma}

\label {l.ct}

\ \\ Set $C = \c_X(T)$. Then

\begin {equation}
C = T \times \Gamma \iso 4^2 \times 3
\end {equation}

where $\Gamma = \la \gamma \ra \iso 3$.

\end {lemma}

\begin {proof}

Let $z$ be one of the three involutions in $T$ and $H = \c_X(z)$. By \zitL, $X$ possesses exactly one conjugacy class of elements of order 2, and their centralizers are isomorphic to $2 \Schur A_{11}$. Thus, as $T$ is abelian, $T \le H \iso 2 \Schur A_{11}$ and $C = \c_H(T)$.

It is easy to see that in $H$ all subgroups isomorphic to $4^2$ are conjugate and have centralizers in $H$ of the form $T \times \Gamma$ with $|\Gamma| = 3$. Obviously $C = C_H(T)$ contains only two elements of order three, each generating $\Gamma$. Direct calculation shows that $\gamma$ centralizes $T$. Therefore $\Gamma = \la \gamma \ra$, and the proof is complete.

\end {proof}

The following result is trivial but will be needed soon.

\begin {lemma}

\label {l.oat}

\ \\ The automorphism group of $T$ has order $96 = 2^5 \cdot 3$.

\end {lemma}

\begin {proof}

$T$ is generated by any pair of elements $t_1, t_2$ of order 4 with $t_1^2 \ne t_2^2$. Thus each automorphism of $T$ must map $(t_1,t_2)$ to an arbitrary pair obeying the same conditions, and every such mapping extends uniquely to an element of $\aut T$. For the image of $t_1$ there are exactly 12 possibilities, and then the image of $t_2$ may be chosen freely among 8 candidates.

\end {proof}

We list a few important informations about $T$, $N$ and their mutual relationship.

\begin {lemma}

\label {l.pi}

\ \\ Let $(P,Q) \in \frak L$ be an oriented line in the apartment $\frak A$. Then

\begin {enumerate}

\item $R(L)$ and $r(L)$ permute the set $\menge {X_{L'}} {L' \in \frak L}$. This extends uniquely to an epimorphism $\pi: N \longmapsto \aut(\frak A) \iso W$.

\item The image of $R(L)$ under $\pi$ interchanges the two numbers which are not among the names of $P$ and $Q$ and the two letters which are not coordinates of $P$.

\item The image of $r(L)$ under $\pi$ interchanges the two numbers not occurring in the names of $P$ and $Q$.

\end {enumerate}

\end {lemma}

\begin {proof}

By Lemma \ref {l.wort} it suffices to verify the formulas by direct calculation for $L = (1a,2b)$ which does not require much effort. The surjectivity of $\pi$ is then clear because the images of all $R(L)$ and $r(L)$, taken together, span $W$.

\end {proof}

\begin {remark}

For example, $\pi$ maps $R(1a,2b)$ to $(3,4)(b,c)$ and $r(1a,2b)$ to $(3,4)$.

\end {remark}

Recall that by definition the line set $\frak L$ consists of certain $5a$-pure groups of order 5, see \zitNM.

\begin {theorem}

\label {s.tn}

{\ }

\begin {enumerate}

\item The elementwise stabilizer of $\frak L$ under conjugation in $X$ equals $T$; the set stabilizer is $\n_X(T) = N$.

\item $N/T \iso W = S_4 \times S_3$.

\item $|N| = 2304 = 2^8 \cdot 3^2$.

\end {enumerate}

\end {theorem}

\begin {proof}

Set as usual $\frak A = \frak A(T)$ and denote the stabilizers of $\frak L$ by $T'$ (elementwise) and $N'$ (setwise). Lemma \ref {l.pi} implies $N \le N'$. Furthermore, by definition $T' \nteq N'$. The quickly verified relations

\begin {align}
t_1(1a,2b) &= R(4a,1b)^{-1} R(1a,4b) \\
t_2(1a,2b) &= R(2a,4b)^{-1} R(4a,2b) \\
t_3(1a,2b) &= R(1a,2b)^{-1} R(2a,1b)
\end {align}

establish that $T = \la t_1(1a,2b), t_2(1a,2b) \ra$ is a subgroup of $N$ and of $T'$. The latter property stems from the fact that the two $R$-factors in each equation correspond to the same involution in $W$.

Every element of $T'$ normalizes each $X_L$, $L \in \frak L$, and induces on it an element of $\aut X_L \iso \aut 5 \iso 4$. This provides us with a homomorphism $\psi$ from $T'$ to the direct product

\begin {equation}
 {\underset {L \in \frak L} {\text {\sffamily \Huge X}}} \aut X_L \iso 4^{36}
\end {equation}

If $\tau$ lies in the kernel of $\psi$, it centralizes all $X_L$ and hence also $\la X_L \ra = X$. This means $\tau \in \c_X(X) = \z(X) = 1$. Thus $\psi$ is injective, and we get $T' \iso \psi(T') \text { \small {$\lesssim$} } 4^{36}$. In particular, $T'$ is a 2-group of order $2^m$, say.

The natural action of $N'$ on $\frak L$ provides us with a homomorphism $\pi: N' \longmapsto \aut \frak A \iso W = S_4 \times S_3$ which is surjective because the subgroup $N$ already gives $\pi(N) = \aut \frak A$. This implies $N'/T' \iso W$.

From these facts we deduce

\begin {equation}
2^{m+4} = 16 |T'| \,\big|\, 144 |T'| = |N'| \,\big|\, |X| = |Ly| = 2^8 \cdot 3^7 \cdot 5^6 \cdot 7 \cdot 11 \cdot 31 \cdot 37 \cdot 67
\end {equation}

leading to $m \le 4$ and $|T'| \le 16$. Since $T \le T'$ has the same order, we arrive at $T' = T$. The restriction $\pi|_N$ has kernel $T' \cap N = T$ and image $W$, whence $|N| = |T| \cdot |W| = 16 \cdot 144 = 2304 = 2^8 \cdot 3^2$.

This implies $T = T' \nteq N'$ and then $T \nteq N$ since $T \le N \le N'$. Both $N$ and $N'$ therefore have the normal subgroup $T$ with quotient $W$ and consequently are of the same order. From $N \le N'$ we obtain $N' = N$.

The normality of $T$ in $N$ is tantamount to $N \le \n_X(T)$. On the other hand, $\c_X(T)$ is also normal in $\n_X(T)$ with factor group isomorphic to some subgroup of $\aut T$. With the help of Lemma \ref {l.ct} and Lemma \ref {l.oat} this yields the estimate

\begin {equation}
|\n_X(t)| = |\c_X(T)| \cdot |\n_X(T):\c_X(T)| \,\big|\, 48 \cdot 96 = 2^9 \cdot 3^2
\end {equation}

Making use of the inclusion of $\n_X(T)$ in $X \iso Ly$, this can be slightly but decisively improved to

\begin {equation}
|\n_X(t)| \,\big|\, 2^8 \cdot 3^2 = |N|
\end {equation}

Combining this with $N \le \n_X(T)$, we obtain $\n_X(T) = N$, concluding the proof.

\end {proof}

Finally, it only remains to determine the root system associated with $T$. In fact, we know it already:

\begin {theorem}

\label {s.wurz}

\ \\ The set of root groups for Kantor's torus $T$ in $X \iso Ly$ (to be understood in the sense of \zitNM) equals $\frak L$.

In other words, all $X_L$ with $L \in \frak L$ are root groups relative to the standard torus, and vice versa.

\end {theorem}

\begin {proof}

Let $L \in \frak L$. By Definition \ref {d.x}, $|X_L| = 5$. Moreover, $L$ is contained (as $L_1$) in the star $\mathcal S(L)$. Thus $X_L = \la \Lambda_1 \ra$ in $\frak S(L)$. An application of the star isomorphism with respect to $L$ yields

\begin {equation}
\Lambda^*_1 = \left| \begin {matrix} 1 & . & . \\ . & 1 & 3 \\ . & . & 1 \end {matrix} \right|
\end {equation}

which is clearly normalized by the diagonal group $T^*$ of $SL_3(5)$. This is tantamount to $X_L^T = X_L$.

Carrying this step out for all $L \in \frak L$ and observing that the $X_L$ are pairwise different, we find all 36 root groups, and the proof is complete.

\end {proof}

\section {Geometrical subspaces}

\label {a.geo}

In the above analysis of (the minimal representation of) $Ly$ it is essential to refer to some eigenbasis of a torus $T$ specified in advance. This alone yields the block structure of the elements of $N$ and the root elements as well as their (5-)logarithms.

Clearly, there is a great many of such bases, and one may ask if they are also equivalent in other respects. This is not the case; we may (and do) apply more restrictions of a geometrical spirit in order to achieve some higher symmetry, especially concerning the structure of the blocks themselves.

For that reason, in the present analysis we made a deliberate choice of the concrete presentation (fulfilling the just mentioned conditions). Since it may not be obvious how this was done, we should give a more detailed explanation.

The first observation is that the group $K$ which is fundamental for our construction, can be used to reduce the number of degrees of freedom in a natural way. It is clear that $K$ permutes the irreducible characters of $T$ and thus the eigenspaces $E_I$ as well. This action has three orbits, given by the type of the character, i.\,e. the maximal degree of the character values as roots of unity.

To be precise, the type of $E_I$ is 1 for the first space, 2 for he next three and 4 for the remaining 12. For any $\tau \in \{1,2,4\}$, we can find subgroups $K_\tau$ of $K$ which act simply transitively of the type-$\tau$-spaces. We employ this symmetry in the following way: For all $\tau$ we select one particular eigenspace, e.\,g. $E_{\tau^2}$, define some basis for the latter and transfer this via $K_\tau$ to get the bases of the other eigenspaces of type $\tau$.

This procedure results in a much more uniform description because the number of distinct blocks occurring in the elements of $N$ and the roots is reduced considerably. But we may achieve even more. By a suitable choice of the basis, we may split each of the 16 natural sections considered above into "minisections" according to the action of $N$. This yields the following refinement:

The section $\sigma_1$ of length 9 decomposes into 4 minisections of lengths $(1^3,6)$ where the exponent describes the multiplicity, the next three section (length 6) into $(2,4)$ and the last 12 sections (length 7) into $(1,2^3)$. Referring to a basis respecting this more detailed decomposition, all elements of $N$ are "miniblock-monomial" with respect to the miniblocks given by the in all 58 minisections with length distribution $(1^{15}, 2^{39},4^3,6)$. However, this does not hold for the root logarithms; they only "see" the eigenspaces, but no finer structure.

Needless to say, the concrete basis of the present paper fulfills all the requirements discussed here. In fact, we even applied several additional criteria, the description of which would lead us too far astray.

In this section, for the convenience of the reader, we want to give a complete enumeration of the minimal $H$-invariant subspaces of $V$, where $H$ is some group of geometrical importance. We explicitly discuss these features for the Sylow-5-group $S$ of $X$, its normalizer (Borel group), the torus normalizer $N$ and the objects of a particular maximal flag in the sense of Kantor.

The constructions are elementary and consist of short and simple direct calculations and repeated applications of the following easy lemma.

\begin {lemma}

\label {l.fix}

\ \\Consider a finite-dimensional vector space $W \ne 0$ over some Galois field $\ds F_q$ of characteristic $p$ (hence $q$ is a power of the prime $p$) and a $p$-group $Q \le GL(W)$.

\begin {enumerate}

\item The fixed space $\Phi = \fix_W(Q)$ of $Q$ in $W$ has positive dimension.

\item Assume in addition that $\dim \Phi = 1$ and $Q \le H \le GL(W)$. Then the closure of all $H$-conjugates of $\Phi$, namely

\begin {equation}
M_H = \gruppe {\fix_W(Q^x)} {x \in H} = \gruppe {\Phi \cdot x} {x \in H}
\end {equation}

is the unique (nonzero) minimal $H$-invariant subspace in $W$.

\item Under the same condition, we furthermore have:

With respect to an arbitrary nonzero $H$-invariant quadratic form on $W$, either $M_H$ is isotropic or $H$ is irreducible.

\end {enumerate}

\end {lemma}

\begin {proof}

\begin {enumerate}

\item The set $W \ohne \Phi$ of nonfixed vectors is permuted by $Q$ in orbits whose lengths are greater than 1 and divisors of $|Q|$, that means $p$-powers. Consequently they are multiples of $p$. The last property is shared by $W$ since $|W| = q^{\dim W} > 1$ likewise is a power of $p$, and hence also by the difference set $\Phi$. Moreover $0 \in \Phi$; so $|\Phi|$ is a positive multiple of $p$ and in particular $|\Phi| \ge p > 1$ which is equivalent to $\dim \Phi > 0$.

\item Applying this to a minimal $H$-invariant subspace $U \ne 0$, we find $1 \le \dim \fix_U(Q) = \dim (U \cap \Phi) \le \dim \Phi = 1$ and further $\dim (U \cap \Phi) = \dim \Phi$ which implies $U \cap \Phi = \Phi$ and $\Phi \subseteq U$. Hence $M_H$, defined via the formula above, is a nonzero $H$-invariant space contained in $U$. By the minimality of $U$ the claim follows.

\item The orthogonal subspace of $U$ in $W$ with respect to the given quadratic form will be written $U^\perp$. Clearly, invariance of $U$ under $H$ implies invariance of $U^\perp$. Isotropy of $U$ means $U \subseteq U^\perp$. The fact that the form does not vanish identically is tantamount to $W^\perp \ne W$, and the irreducibility of $H$ to $M_H = W$. Thus we have to show:

{\it Provided $W^\perp \ne W$, either $M_H = W$ or $M_H \subseteq M_H^\perp$.}

In order to prove this, assume $W^\perp \ne W$. If $M_H^\perp \ne 0$, then $M_H^\perp$ is $H$-invariant and must contain a minimal space with the same property, which can only be $M_H$; hence $M_H \subseteq M_H^\perp$.

On the other hand, $M_H^\perp = 0$ implies $W^\perp \subseteq M_H^\perp = 0$ and then $W^\perp = 0$. This is equivalent to the nondegeneracy of the given quadratic form, and we get $M_H = M_H^{\perp\perp} = 0^\perp = W$.

Hence (at least) one of the two possibilities always holds true. To see that both are incompatible with each other, we note that $M_H = W$ and $M_H \subseteq M_H^\perp$ together first yield $M_H^\perp = W$ and then $W^\perp = W$, contradicting the hypothesis.

\end {enumerate}

\end {proof}

For the sake of clarity we next introduce a few more pieces of notation:

\begin {definition}

\label {d.ech}

{\ }

\begin {enumerate}

\item The {\bf echelon matrix} of a subspace $U \subseteq V$ is denoted by $\frak B(U)$. Occasionally, if we need not calculate with the matrix, we simply write $\frak B(U)$ as list of its row vectors.

\item We call a subspace $U \subseteq V$ {\bf special} if it possesses a basis of unit vectors. In that case, the echelon basis $\frak B(U)$ has the same property, and we use the abbreviation $U = \mathcal S(i_1, \ldots, i_d)$ for the space spanned by $\mng {e_{i_1}, \ldots, e_{i_d}}$. More precisely, the echelon basis of $U$ is

\begin {equation}
\frak B(U) = \vkt {e_{i_1}, \ldots, e_{i_d}}
\end {equation}

where the subscripts have to be arranged in natural (ascending) order.

\item For $1 \le I \le 16$, the echelon basis of $E_I$ is $\vektor {e_i} {\alpha_I \le i \le \omega_I}$. Thus each $u \in E_I$ has a unique representation of the shape

\begin {equation}
u = (u_{\alpha_I} \ldots u_{\omega_I}) \cdot \frak B(E_I) = (0 \ldots 0 | u_{\alpha_I} \ldots u_{\omega_I} | 0 \ldots 0)
\end {equation}

which will simply be written in shorthand notation as

\begin {equation}
u = [u_{\alpha_I} \ldots u_{\omega_I}]_I
\end {equation}

\end {enumerate}

\end {definition}

We want to give a few examples:

\begin {enumerate}

\item By direct calculation one finds that the fixed space of $S \in \syl_5(X)$ is one-dimensional and contained in $E_5$. It has the echelon matrix

\begin {equation}
\frak B(\fix S) = \left( [1 \ . \ 3 \ . \ . \ 4 \ 2]_5 \right)
\end {equation}

\item The eigenspaces $E_I$, $1 \le I \le 16$, of our standard torus $T$ are obviously special:

\begin {equation}
E_I = \mathcal S(\alpha_I, \alpha_I+1, \ldots, \omega_I-1, \omega_I)
\end {equation}

Several other special subspaces of $V$ will be introduced below.

\end {enumerate}

\begin {theorem}

\label {s.mini}

{\ }

\begin {enumerate}

\item Consider the maximal flag $(P,L,F)$ in $\frak A(T)$ with $P = (1a) \in \frak P$, $L = (1a,2b) \in \frak L$ and $F = (1a,2b,3c) \in \frak F$. The intersection of the associated Kantor objects $O_K = \n_X(O)$, $O \in \{P,L,F\}$, coincides with the Borel group

\begin {equation}
P_K \cap L_K \cap F_K = \n_X(S) = S:T
\end {equation}

where $S$ is the unique Sylow-5-subgroup of $X$ shared by all flag components (in the sense of Kantor).

Every $H \in \{S,{S:T},P_K,L_K,F_K\}$ possesses a unique minimal invariant subspace $M_H \ne 0$ in $V$. The echelon bases of these $M_H$ are

\begin {equation}
\frak B(M_S) = \frak B(M_{S:T}) = \big( [1 \ . \ 3 \ . \ . \ 4 \ 2]_5 \big)
\end {equation}

\hoch

\begin {align}
\frak B(M_{P_K}) = \big( &[. \ . \ 1 \ 2 \ 4 \ . \ . \ . \ .]_1, [1 \ . \ 3 \ . \ . \ 4 \ 2]_5, [1 \ . \ 3 \ . \ . \ 1 \ 3]_7, [1 \ . \ 3 \ . \ . \ 4 \ 2]_9, \nonumber \\
&[1 \ . \ 3 \ . \ . \ 1 \ 3]_{11}, [1 \ . \ 3 \ . \ . \ 4 \ 2]_{13}, [1 \ . \ 3 \ . \ . \ 1 \ 3]_{15} \big)
\end {align}

\hoch

\begin {equation}
\frak B(M_{L_K}) = \big( [1 \ . \ 3 \ . \ . \ 4 \ 2]_5, [1 \ 2 \ 2 \ . \ . \ 3 \ 2]_{12}, [1 \ 2 \ 2 \ . \ . \ 2 \ 3]_{14}, [1 \ . \ 3 \ . \ . \ 1 \ 3]_{15} \big)
\end {equation}

\hoch

\begin {align}
\frak B(M_{F_K}) = \big( &[. \ 1 \ . \ 2 \ . \ . \ 2 \ 3 \ 3]_1, [1 \ . \ 3 \ . \ . \ 4 \ 2]_5, [1 \ . \ 3 \ . \ . \ 1 \ 3]_7, [1 \ 3 \ . \ . \ . \ 2 \ 4]_8, [1 \ 3 \ . \ . \ . \ 2 \ 4]_9, \nonumber \\
&[1 \ 2 \ 2 \ . \ . \ 3 \ 2]_{12}, [1 \ 2 \ 2 \ . \ . \ 2 \ 3]_{13}, [1 \ 2 \ 2 \ . \ . \ 2 \ 3]_{14}, [1 \ . \ 3 \ . \ . \ 1 \ 3]_{15}, [1 \ 3 \ . \ . \ . \ 3 \ 1]_{16} \big)
\end {align}

\item The special spaces $U_1, \ldots, U_{10}$, defined as $U_m = \mathcal S(J_m)$ for $1 \le m \le 10$ with index sets

\begin {align}
J_1 &= \{ 1 \} \\
J_2 &= \{ 2 \} \\
J_3 &= \{ 3 \} \\
J_4 &= \{ 4, 5, 6, 7, 8, 9 \} \\
J_5 &= \{ 10, 11, 16, 17, 22, 23 \} \\
J_6 &= \{ 12, 13, 14, 15, 18, 19, 20, 21, 24, 25, 26, 27 \} \\
J_7 &= \{ 28, 35, 42, 49, 56, 63, 70, 77, 84, 91, 98, 105 \} \\
J_8 &= (J_7 + 1) \cup (J_7 + 2) \\
J_9 &= (J_7 + 3) \cup (J_7 + 4) \\
J_{10} &= (J_7 + 5) \cup (J_7 + 6)
\end {align}

are invariant under $N$. Relative to the non-degenerate quadratic form $\mathcal F$, the representation module $V$ decomposes into the direct orthogonal sum of the $U_j$.

\item $N$ has exactly 14 minimal invariant spaces, namely $U_1, \ldots, U_{10}$ and the four subspaces $U_{10+m}$ of $U_9 \oplus U_{10}$ with echelon bases

\begin {eqnarray}
\frak E(U_9) + m \frak E(U_{10}) = \liste {e_i + m e_{i+2}} {i \in J_9}
\end {eqnarray}

for $m \in \{1,2,3,4\}$.

\end {enumerate}

\end {theorem}

\begin {proof}

{\ }

\begin {enumerate}

\item In each case $S$ is a 5-subgroup of $H$ with $\dim \fix S = 1$, and Lemma \ref {l.fix} yields the claim.

\item Since concatenation of the echelon bases of $U_1, \ldots, U_{10}$ results in a row permutation of the unit matrix, we get $V = U_1 \oplus \ldots \oplus U_{10}$. A simple calculation furthermore shows that all $U_m$ are invariant under the elements $n_1, \ldots, n_4$ of Lemma \ref {l.Nerz} and thus also under $\la n_1, n_2, n_3, n_4 \ra = N$.

Finally, let $I,J \in \{ 1, \ldots, 10 \}$ with $I \ne J$ and $i,j \in \{ 1, \ldots, 111 \}$ such that $e_i \in U_I$ and $e_j \in U_J$. Then

\begin {equation}
\la e_i,e_j \ra = e_i \cdot f \cdot \tr e_j = f_{ij} = 0
\end {equation}

and hence

\begin {equation}
\frak E(U_I) \cdot f \cdot \tr \frak E(U_J) = 0
\end {equation}

which amounts to the same thing as $U_I \bot U_J$.

\item In order to prove the irreducibility of the $U_j$, $1 \le j \le 10$, as $N$-modules, we use character theory. First we need the character table of $N$. Unfortunately, the elements of $N$ are (111,111)-matrices and therefore too clumsy for a direct application of the Dixon-Schneider algorithm which is part of the {\sc Gap}-library.

But a short detour helps. As we know already, $N$ permutes the 144 root elements (transitively). This provides us with a homomorphism $\kappa: N \longrightarrow S_{144}$.

Any element of the kernel of $\kappa$ centralizes all root elements and thus the group $X$ generated by them. Consequently, the kernel is contained in $\c_X(X) = \z(X) = 1$, and $\kappa$ is injective. Hence $N \iso \kappa(N)$, and the characters of the latter group are easily and quickly found by the proposed method.

From the character tables of $N$ and $Ly$ (see \zitL), the fusion map between both groups is essentially unique and can be calculated without difficulty, even by hand. This also gives the restriction $\psi = \phi|_N$, which is the character of the $N$-action on $V$.

These results are shown in Tab. \ref {t.char}. The first row contains the {\sc Atlas} names of the 30 conjugacy classes of $N$; underneath (second row) the corresponding $Ly$-classes are listed. The next 30 lines consist of the irreducible characters of $N$, while the last row comprises the values of $\psi$.

Since $N$ (order 2304) is a $5'$-group, the representation of $N$ and its Brauer-5-character $\psi$ may be interpreted as ordinary ones. In particular, Maschke's theorem holds and Frobenius theory yields the decomposition of $\psi$ into irreducible constituents. More precisely, $\psi = \psi_1 + \ldots + \psi_{10}$ with degrees $111 = 1 + 1 + 1 + 6 + 6 + 12 + 12 + 24 + 24 + 24$.

All $\psi_l$ are distinct, except for two characters of degree 24, which we may identify with $\psi_9$ and $\psi_{10}$.

\begin {center}

\begin {tabelle}

\label {t.char}

\begin {center}
\hoch
Character table of $N = \n_X(T) = T \cdot W$ \hoch \
\end {center}

{\fontsize {10} {12} \selectfont \setlength {\tabcolsep} {2.2pt}
\begin {tabular*} {17cm} {*{30}{r}} \\ \hline
1a & 2a & 2b & 2c & 2d & 2e & 3a & 3b & 3c & 4a & 4b & 4c & 4d & 4e & 6a & 6b & 6c & 6d & 6e & 8a & 8b & 8c & 8d & 8e & 12a & 12b & 12c & 12d & 12e & 24a \\
1a & 2a & 2a & 2a & 2a & 2a & 3a & 3b & 3b & 4a & 4a & 4a & 4a & 4a & 6a & 6a & 6a & 6a & 6c & 8a & 8a & 8b & 8b & 8b & 12a & 12a & 12a & 12a & 12a & 24a \\ \hline
1 & 1 & 1 & 1 & 1 & 1 & 1 & 1 & 1 & 1 & 1 & 1 & 1 & 1 & 1 & 1 & 1 & 1 & 1 & 1 & 1 & 1 & 1 & 1 & 1 & 1 & 1 & 1 & 1 & 1 \\
1 & 1 & 1 & -1 & -1 & 1 & 1 & 1 & 1 & 1 & 1 & -1 & 1 & -1 & 1 & 1 & 1 & -1 & -1 & -1 & -1 & -1 & 1 & 1 & 1 & 1 & -1 & -1 & -1 & -1 \\
1 & 1 & 1 & -1 & 1 & -1 & 1 & 1 & 1 & 1 & 1 & -1 & -1 & -1 & 1 & 1 & 1 & -1 & 1 & -1 & 1 & 1 & -1 & -1 & 1 & 1 & -1 & -1 & -1 & -1 \\
1 & 1 & 1 & 1 & -1 & -1 & 1 & 1 & 1 & 1 & 1 & 1 & -1 & 1 & 1 & 1 & 1 & 1 & -1 & 1 & -1 & -1 & -1 & -1 & 1 & 1 & 1 & 1 & 1 & 1 \\
2 & 2 & 2 & -2 & . & . & -1 & 2 & -1 & 2 & 2 & -2 & . & -2 & -1 & -1 & -1 & 1 & . & -2 & . & . & . & . & -1 & -1 & 1 & 1 & 1 & 1 \\
2 & 2 & 2 & 2 & . & . & -1 & 2 & -1 & 2 & 2 & 2 & . & 2 & -1 & -1 & -1 & -1 & . & 2 & . & . & . & . & -1 & -1 & -1 & -1 & -1 & -1 \\
2 & 2 & 2 & . & -2 & . & 2 & -1 & -1 & 2 & 2 & . & . & . & 2 & 2 & 2 & . & 1 & . & -2 & -2 & . & . & 2 & 2 & . & . & . & . \\
2 & 2 & 2 & . & 2 & . & 2 & -1 & -1 & 2 & 2 & . & . & . & 2 & 2 & 2 & . & -1 & . & 2 & 2 & . & . & 2 & 2 & . & . & . & . \\
3 & 3 & -1 & -1 & 3 & -1 & 3 & . & . & 3 & -1 & -1 & -1 & 1 & 3 & -1 & -1 & -1 & . & -1 & -1 & -1 & -1 & 1 & 3 & -1 & -1 & 1 & 1 & -1 \\
3 & 3 & -1 & -1 & -3 & 1 & 3 & . & . & 3 & -1 & -1 & 1 & 1 & 3 & -1 & -1 & -1 & . & -1 & 1 & 1 & 1 & -1 & 3 & -1 & -1 & 1 & 1 & -1 \\
3 & 3 & -1 & 1 & 3 & 1 & 3 & . & . & 3 & -1 & 1 & 1 & -1 & 3 & -1 & -1 & 1 & . & 1 & -1 & -1 & 1 & -1 & 3 & -1 & 1 & -1 & -1 & 1 \\
3 & 3 & -1 & 1 & -3 & -1 & 3 & . & . & 3 & -1 & 1 & -1 & -1 & 3 & -1 & -1 & 1 & . & 1 & 1 & 1 & -1 & 1 & 3 & -1 & 1 & -1 & -1 & 1 \\
4 & 4 & 4 & . & . & . & -2 & -2 & 1 & 4 & 4 & . & . & . & -2 & -2 & -2 & . & . & . & . & . & . & . & -2 & -2 & . & . & . & . \\
6 & 6 & -2 & -2 & . & . & -3 & . & . & 6 & -2 & -2 & . & 2 & -3 & 1 & 1 & 1 & . & -2 & . & . & . & . & -3 & 1 & 1 & -1 & -1 & 1 \\
6 & 6 & -2 & 2 & . & . & -3 & . & . & 6 & -2 & 2 & . & -2 & -3 & 1 & 1 & -1 & . & 2 & . & . & . & . & -3 & 1 & -1 & 1 & 1 & -1 \\
6 & 6 & 2 & -2 & . & . & 6 & . & . & -2 & -2 & -2 & . & . & 6 & 2 & 2 & -2 & . & 2 & . & . & . & . & -2 & -2 & -2 & . & . & 2 \\
6 & 6 & 2 & 2 & . & . & 6 & . & . & -2 & -2 & 2 & . & . & 6 & 2 & 2 & 2 & . & -2 & . & . & . & . & -2 & -2 & 2 & . & . & -2 \\
6 & 6 & -2 & . & . & -2 & 6 & . & . & -2 & 2 & . & -2 & . & 6 & -2 & -2 & . & . & . & . & . & 2 & . & -2 & 2 & . & . & . & . \\
6 & 6 & -2 & . & . & 2 & 6 & . & . & -2 & 2 & . & 2 & . & 6 & -2 & -2 & . & . & . & . & . & -2 & . & -2 & 2 & . & . & . & . \\
6 & 6 & 2 & -2 & . & . & -3 & . & . & -2 & -2 & -2 & . & . & -3 & A & *A & 1 & . & 2 & . & . & . & . & 1 & 1 & 1 & B & -B & -1 \\
6 & 6 & 2 & -2 & . & . & -3 & . & . & -2 & -2 & -2 & . & . & -3 & *A & A & 1 & . & 2 & . & . & . & . & 1 & 1 & 1 & -B & B & -1 \\
6 & 6 & 2 & 2 & . & . & -3 & . & . & -2 & -2 & 2 & . & . & -3 & A & *A & -1 & . & -2 & . & . & . & . & 1 & 1 & -1 & -B & B & 1 \\
6 & 6 & 2 & 2 & . & . & -3 & . & . & -2 & -2 & 2 & . & . & -3 & *A & A & -1 & . & -2 & . & . & . & . & 1 & 1 & -1 & B & -B & 1 \\
12 & 12 & -4 & . & . & . & -6 & . & . & -4 & 4 & . & . & . & -6 & 2 & 2 & . & . & . & . & . & . & . & 2 & -2 & . & . & . & . \\
12 & -4 & . & 2 & . & -2 & 12 & . & . & . & . & -2 & 2 & . & -4 & . & . & 2 & . & . & -2 & 2 & . & . & . & . & -2 & . & . & . \\
12 & -4 & . & 2 & . & 2 & 12 & . & . & . & . & -2 & -2 & . & -4 & . & . & 2 & . & . & 2 & -2 & . & . & . & . & -2 & . & . & . \\
12 & -4 & . & -2 & . & 2 & 12 & . & . & . & . & 2 & -2 & . & -4 & . & . & -2 & . & . & -2 & 2 & . & . & . & . & 2 & . & . & . \\
12 & -4 & . & -2 & . & -2 & 12 & . & . & . & . & 2 & 2 & . & -4 & . & . & -2 & . & . & 2 & -2 & . & . & . & . & 2 & . & . & . \\
24 & -8 & . & 4 & . & . & -12 & . & . & . & . & -4 & . & . & 4 & . & . & -2 & . & . & . & . & . & . & . & . & 2 & . & . & . \\
24 & -8 & . & -4 & . & . & -12 & . & . & . & . & 4 & . & . & 4 & . & . & 2 & . & . & . & . & . & . & . & . & -2 & . & . & . \\ \hline
111 & -1 & -1 & -1 & -1 & -1 & -24 & 3 & 3 & 3 & 3 & 3 & 3 & 3 & 8 & 8 & 8 & 8 & -1 & -3 & -3 & 1 & 1 & 1 & . & . & . & . & . & . \\ \hline
\end {tabular*}
} \hoch

\begin {align*}
A &= - 1 - 2 \sqrt{-3} \\
*A &= - 1 + 2 \sqrt{-3} \\
B &= \sqrt{-3}
\end {align*}

\end {tabelle}

\end {center}

The coincidence of the character degrees with the dimensions of $U_1, \ldots, U_{10}$ proves that the direct sum decomposition given above cannot be refined and the $U_j$ are irreducible as $N$-spaces.

Furthermore, every minimal $N$-invariant subspace of $V$ other than the $U_j$ must be contained in $U_9 \oplus U_{10}$ and isomorphic as $N$-module to $U_9 \iso U_{10}$. This immediately leads to the proposition, at the same time concluding the proof of the theorem.

\end {enumerate}

\end {proof}

Let $z \in 2a$ be any involution in $X$. Then the eigenspaces

\begin {equation}
V_\pm(z) = \menge {u \in V} {u z = \pm u}
\end {equation}

of $z$ with eigenvalues $\pm 1$ are clearly invariant under $B = \c_X(z) = \n_X(z) \iso 2 \Schur A_{11}$. The fixed space $V_+(z)$ is 55-dimensional, the minus space $V_-(z)$ has dimension 56. They are orthogonal complements of each other.

By inspection of the 5-Brauer characters of $2 \Schur A_{11}$, cf. \zitMN, both are irreducibly acted upon by $B$ and thus the only nontrivial $B$-invariant subspaces of $V$.

If $z \in T$, that means $z = t_i(1a,2b)^2$ for some $i \in \{1,2,3\}$, these spaces are special, namely of the form $\mathcal S(I_\pm)$ with index sets defined as

\begin {equation}
I_\pm = \menge {i \in \{1, \ldots, 111\}} {e_i \in V_\pm(z)} = \menge {i \in \{1, \ldots, 111\}} {e_i z = \pm e_i}
\end {equation}

Among the two $Ly$-classes $3a$ and $3b$ of order-3-elements, the former is by far the more interesting and geometrically important, and we restrict our attention to it.

Consider $y \in 3a$. Then $\c_X(y) \iso 3 \Schur Mc$ and $\n_X(y) \iso (3 \Schur Mc):2$. Evidently the 21.dimensional fixed space $\fix(y)$ of $y$ and its orthogonal complement $\fix(y)^\bot$ of dimension 90 are invariant under both groups.

The restriction of $\phi$ to the normalizer of $y$ decomposes (see {\it loc.\,cit.}) into two irreducible components of degrees 21 and 90, respectively.

Reducing further to $\c_X(y)$, the fixed space of $y$ remains irreducible, while the 90-dimensional character splits into two irreducible constituents of degree 45 each forming a complex conjugate pair.

Consequently, after some sufficiently large extension of the underlying field we get two nontrivial invariant 45-dimensional invariant subspaces: the eigenspaces of $y$ whose eigenvalues are the proper third roots of unity. Since we are calculating over the ground field $\F_5$ itself, however, in $V$ this cannot happen (the smallest splitting field would be $\F_{25}$).

In conclusion: Under both groups, $\c_X(y) \iso 3 \Schur Mc$ as well as $\n_X(y) \iso (3 \Schur Mc):2$, there exist only two nontrivial invariant subspaces of $V$, namely $\fix(y) \iso \F_5^{21}$ and $\fix(y)^\bot \iso \F_5^{90}$. Each is the orthogonal complement of the other.

With respect to our geometry, the particular example $y \in \c_X(T) = T \times \Gamma$, that means $y = \gamma^\pm$, is of major interest. In that case both spaces are special. The fixed space of $\gamma$ (or $\gamma^{-1}$) is generated by

\begin {equation}
\menge {e_i} {i \in \{  1, 2, 3, 10, 11, 16, 17, 22, 23, 28, 35, 42, 49, 56, 63, 70, 77, 84, 91, 98, 105 \}}
\end {equation}

and its orthogonal complement by the unit vectors not contained in this list.

These few examples abundantly demonstrate that and how our deliberate choice of the underlying basis of the module $V \iso \F_5^{111}$ helps to clarify the structure of the Lyons group and its minimal 5-representation as well as the associated Kantor geometry. We confidently hope that in the long run our approach will lead to a much better understanding of these entities.

Perhaps similar Ans\"atze might also be successful for (the) other sporadic groups. Time will tell \ldots

\section {Summary and outlook}

\label {a.summ}

During the last three or four decades, those 20 out of the 26 sporadic groups which belong to the Mathieu-Conway-Fischer-Monster family have been the subject of a great many of thorough investigations from several points of view, while the remaining 6 "exotic" groups did not attract much interest, but were treated as Cinderella subjects. This, of course, is not justified since for gaining a deeper insight into the phenomenon "simple group" all sporadics are equally important.

The aim of the present investigation therefore is to offer a modest contribution to a more thorough understanding of the geometry of the second largest exotic, namely the Lyons group $Ly$ which among all sporadic groups seemingly is by far the one most closely related to the Chevalley groups.

The original proofs of existence and uniqueness of the group itself as well as its absolute minimal representation (111-dimensional over $\F_5$) use generators (different for both problems) which are in no way adapted to the beautiful geometric properties of $Ly$, and the same holds for the defining relations which are chosen in an {\it ad hoc} manner. This causes many inconveniences and is largely responsible for the length of the rather cubersome methods to be applied.

In particular, the Sims-Havas presentation with 5 generators and 86 more or less longish and unnatural relations requires heavy computer calculations (Todd-Schreier algorithm) which were at the border of the technical possibilities in the 1960-ies. In fact, it turned out necessary to invent some new ideas in order to overcome the difficulties. The same must be said concerning the construction of the representation.

Hence the major goal of this paper is to unify the approach by constructing $Ly$ as a subgroup $X$ of $SL_{111}(5)$. This simultaneously demonstrates the existence of the group. Furthermore we employ the marvelous properties of the Kantor geometry in its extended form by choosing a root system (i.\,e. a special set of 36 elements, one in each root subgroup) as generators. This provides us with a completely symmetric presentation, the relations of which are given by the Chevalley-Steinberg theory, thus minimizing the arbitraryness of the construction as much as possible.

The larger number of generators as compared to Sims' approach is outweighed by the symmetry, and we make use of a recent characterization of the Lyons group due to Gr\"uninger. This enables us to avoid the cumbersome coset enumeration algorithm in favour of a much more natural amalgam construction of the Lyons group. Moreover, the difficult search for specific subgroups needed to establish the representation becomes completely obsolete.

The Ansatz followed in the present investigation is to choose a Kantor torus $T$ and express all group elements with respect to an eigenbasis of $T$. This leads to a decomposition of the matrices into $16 \times 16$ blocks which are easy to handle.

The normalization just described makes all elements of the torus normalizer $N = \n_X(T)$ block-monomial. This by itself is not remarkable since $N$ is a tiny subgroup in $X \iso Ly$. Nevertheless, the choice of a $T$-eigenbasis for the representation brings an additional bonus. It came as a surprise (and I do not have a natural explanation for it) that relative to an arbitrary eigenbasis the logarithms of the roots (taken modulo 5) are also block-monomial.

The latter property allows to construct $X$ (or $Ly$) with the help of only four block-monomial matrices; a fifth is given in addition which is not reqired for the group construction, but defines a nondegenerate $X$-invariant quadratic form which is useful in many respects.

By some experimentation with {\sc Gap} a particular basis was found such that the minimal subspaces invariant under several geometrically interesting and thus important subgroups become most convenient.

It should be noted that the approach given here not only avoids methods of combinatorial group theory completely; it reduces the amount of computational work to a minimum. In fact, the application of computers is restricted to elementary matrix algebra, i.\,e. multiplication and addition of 111-dimensional matrices over the ground field with 5 elements.

It is to be hoped that the new representation and its proof of correctness will stimulate further research into the subject and lead to more knowledge of the extremely fascinating Lyons group (and perhaps the other 5 exotics as well).

\section {Acknowledgements}

\label {a.ack}

All explicit calculations within the Lyons group and some related objects have been carried out with the help of the computer algebra package {\sc Gap} invented by Joachim Neub\"user and developed to perfection over the years by himself, his colleagues, assistants, students and external collaborators.

It is not only my obligation, but also a great pleasure to acknowledge my gratitude to all {\sc Gap} developers at RWTH Aachen for invaluable support in many respects. Over more than three decades, I was always welcomed as a friend, and even inconvenient or stupid questions were answered with great enthusiasm.

When I had problems concerning calculations within finite groups which I could not solve at Bonn, the standard procedure "go to Aachen and ask some expert" quite often helped to cut the Gordian knot.

Beyond this, I highly appreciate the personal relations with several members of the {\sc Gap} developers' community.

I am most grateful to my former student Ulrich Kaiser (Aachen) who implemented the package under Windows XP, thus enabling me to do research conveniently at home. This was an indispensable help to accelerate the investigations.

The manuscript of this publication was written with LaTeX on a laptop. For the sake of efficiency this required to install {\sc Gap} on the basis of the operating system Linux.

I am deeply indebted to Thomas Breuer (RWTH) who spent much time and energy on a "remote diagnosis" of my computing device, eventually eliminating several bugs, among them some serious ones.

So, my friends, thanks to all of you\,! This research is also yours; I hope the result will please you and justify the efforts.

\section {References}

\label {a.ref}

Roger W. {\bf Carter} {\bf [1972]} \\
{\it Simple Groups of Lie type} \\
Wiley \& Sons, London / New York / Sydney / Toronto \\

John H. {\bf Conway} , Robert T. {\bf Curtis}, Simon P. {\bf Norton}, Richard A. {\bf Parker}, Robert A. {\bf Wilson} {\bf [1984]} \\
{\it An Atlas of finite groups} \\
Clarendon Press, Oxford \\

Matthias {\bf Gr\"uninger} {\bf [2007]} \\
{\it Die Geometrie der Lyonsgruppe} \\
Dissertation T\"ubingen \\

George {\bf Havas}, Charles C. {\bf Sims} {\bf [1999]} \\
{\it A presentation for the Lyons simple group} \\
{\it in} {P. Dr\"axler, C.\,M. Ringel, G.\,O. Michler (eds.): {\it Computational methods for representations of groups and algebras}, 241--249 \\
Birkh\"auser, Basel \\

James E. {\bf Humphreys} {\bf [1975]} \\
{\it Linear algebraic groups} \\
(Graduate Texts in Mathematics, vol. \ul {21}) \\
Springer Verlag, Berlin, Heidelberg, New York \\

Christopher {\bf Jansen}, Robert A. {\bf Wilson} {\bf [1996]} \\
{\it The minimal faithful 3-modular representation for the Lyons group} \\
Comm. Alg. \ul {24}, 873--879 \\

William {\bf Kantor} {\bf [1981]} \\
{\it Some geometries that are almost buildings} \\
Eur. J. Comb. \ul {2}, 239--247 \\

Richard {\bf Lyons} {\bf [1972]} \\
{\it Evidence for a new finite simple group} \\
J. Alg. \ul {20}, 540--569 \\

Werner {\bf Meyer}, Wolfram {\bf Neutsch} {\bf [1984]} \\
{\it \"Uber 5-Darstellungen der Lyonsgruppe} \\
Math. Ann. \ul {267}, 519--535 \\

Werner {\bf Meyer}, Wolfram {\bf Neutsch}, Richard {\bf Parker} {\bf [1985]} \\
{\it The minimal 5-representation of Lyons' sporadic group} \\
Math. Ann. \ul {272}, 29--39 \\

Eliakim Hastings {\bf Moore} {\bf [1897]} \\
{\it Concerning the abstract groups of order $k!$ and $\tfrac 12 k!$ holohedrically isomorphic with the symmetric and the alternating substitution groups on $k$ letters} \\
Proc. London Math. Soc. \ul {28}, 357--366 \\

Wolfram {\bf Neutsch}, Werner {\bf Meyer} {\bf [1989]} \\
{\it A root system for the Lyons group} \\
Math. Ann. \ul {283}, 285--299 \\

Charles C. {\bf Sims} {\bf [1973]} \\
{\it The existence and uniqueness of Lyons' group} \\
{\it in} T. Gagen, M. Hale, E. Shult (eds.): {\it Finite groups '72}, 138--141 \\
North-Holland, Amsterdam \\

Robert {\bf Steinberg} {\bf [1962]} \\
{\it Generators for simple groups} \\
Canad. J. Math. \ul {14}, 277-283 \\

Robert A. {\bf Wilson} {\bf [1984]} \\
{\it The subgroup structure of the Lyons group} \\
Math, Proc. Cambr. Philos. Soc. \ul {95}, 403--409 \\

Robert A. {\bf Wilson} {\bf [1985]} \\
{\it The maximal subgroups of the Lyons group} \\
Math, Proc. Cambr. Philos. Soc. \ul {97}, 433--436 \\

Andrew J. {\bf Woldar} {\bf [1984]} \\
{\it On the maximal subgroups of Lyons' group and evidence for the existence of a 111-dimensional faithful LyS-module over a field of characteristic 5} \\
Dissertation Ohio \\

Andrew J. {\bf Woldar} {\bf [1987]} \\
{\it On the maximal subgroups of Lyons' group} \\
Comm. Alg. \ul {15}, 1195--1203

\end {document}